\documentclass[letter,twocolumn]{autart}

\usepackage{amsmath, amssymb, mathtools}
\usepackage{graphicx,color,subfig}
\usepackage{natbib} 
\usepackage{url} 

\newtheorem{theorem}{Theorem}[section]
\newtheorem{proposition}[theorem]{Proposition}

\newtheorem{corollary}[theorem]{Corollary}
\newtheorem{remark}[theorem]{Remark}
\newtheorem{definition}[theorem]{Definition}
\newtheorem{problem}{Problem}

\newcommand{\longthmtitle}[1]{{\rm \textbf{(#1)}}}

\newcommand{\real}{\ensuremath{\mathbb{R}}}
\newcommand{\realpos}{\ensuremath{\mathbb{R}_{>0}}}
\newcommand{\realnonneg}{\ensuremath{\mathbb{R}_{\ge 0}}}
\newcommand{\integerpos}{{\mathbb{N}}}
\newcommand{\integernonneg}{{\mathbb{Z}}_{\ge 0}}

\let\emptyset\varnothing

\newcommand{\Ac}{\mathcal{A}}
\newcommand{\Bc}{\mathcal{B}}

\newcommand{\Gc}{\mathcal{G}}
\newcommand{\Ic}{\mathcal{I}}
\newcommand{\Kc}{\mathcal{K}}
\newcommand{\Lc}{\mathcal{L}}
\newcommand{\Ocal}{\mathcal{O}}
\newcommand{\Pc}{\mathcal{P}}

\newcommand{\Nc}{\mathcal{N}}

\newcommand{\Lap}{\text{Lap}}
\newcommand{\ones}{\mathbf{1}}
\newcommand{\pr}{\mathbb{P}}
\newcommand{\E}{\mathbb{E}}

\newcommand{\tr}{\text{tr}}

\newcommand{\until}[1]{\{1,\dots,#1\}}
\newcommand{\Ave}[1]{\operatorname{Ave}(#1)}

\newcommand{\oprocendsymbol}{\hbox{$\bullet$}}
\newcommand{\oprocend}{\relax\ifmmode\else\unskip\hfill\fi\oprocendsymbol}

\renewcommand{\theenumi}{(\roman{enumi})}

\graphicspath{{figures/}}
\allowdisplaybreaks

\begin{document}

\begin{frontmatter}

 \title{Differentially Private Average Consensus: Obstructions,
    Trade-Offs, and Optimal Algorithm Design}

  \thanks{A preliminary version of this paper appeared
    as~\citep{EN-PT-JC:15-necsys} at the 5th IFAC Workshop on
    Distributed Estimation and Control in Networked Systems.}
    
    \vspace{-10pt}
  
  \author[First]{Erfan Nozari}%
  \author[First]{\quad Pavankumar Tallapragada}%
  \author[First]{\quad Jorge Cort\'es}
  
  \address[First]{Department of Mechanical and Aerospace Engineering,
    University of California, San Diego,
    \{enozari,ptallapragada,cortes\}@ucsd.edu}
  
  \begin{abstract}
    This paper studies the multi-agent average consensus problem under
    the requirement of differential privacy of the agents' initial
    states against an adversary that has access to all the messages.
    We first establish that a differentially private consensus
    algorithm cannot guarantee convergence of the agents' states to
    the exact average in distribution, which in turn implies the same
    impossibility for other stronger notions of convergence.  This
    result motivates our design of a novel differentially private
    Laplacian consensus algorithm in which agents linearly perturb
    their state-transition and message-generating functions with
    exponentially decaying Laplace noise.  We prove that our algorithm
    converges almost surely to an unbiased estimate of the average of
    agents' initial states, compute the exponential mean-square rate
    of convergence, and formally characterize its differential privacy
    properties. We show that the optimal choice of our design
    parameters (with respect to the variance of the convergence point
    around the exact average) corresponds to a one-shot perturbation
    of initial states and compare our design with various counterparts
    from the literature.  Simulations illustrate our results.
  \end{abstract}
  
  \begin{keyword}
    Average consensus, Differential privacy, Multi-agent systems,
    Exponential mean-square convergence rate, Networked control
    systems
  \end{keyword}
  
\end{frontmatter}

\section{Introduction}

The social adoption of new technologies in networked cyberphysical
systems relies heavily on the privacy preservation of individual
information.  Social networking, the power grid, and smart
transportation are only but a few examples of domains in need of
privacy-aware design of control and coordination strategies.  In these
scenarios, the ability of a networked system to fuse information,
compute common estimates of unknown quantities, and agree on a common
view of the world is critical. Motivated by these observations, this
paper studies the multi-agent average consensus problem, where a group
of agents seek to agree on the average of their individual values by
only interchanging information with their neighbors.  This problem has
numerous applications in synchronization, network management, and
distributed control/computation/optimization. In the context of
privacy preservation, the notion of differential privacy has gained
significant popularity due to its rigorous formulation and proven
security properties, including resilience to post-processing and side
information, and independence from the model of the adversary.
Roughly speaking, a strategy is differentially private if the
information of an agent has no significant effect on the aggregate
output of the algorithm, and hence its data cannot be inferred by an
adversary from its execution.  This paper is a contribution to the
emerging body of research that studies privacy preservation in
cooperative network systems, specifically focused on gaining insight
into the achievable trade-offs between privacy and performance in
multi-agent average consensus.

\emph{Literature Review:} The problem of multi-agent average consensus
has been a subject of extensive research in networked systems and it
is impossible to survey here the vast amount of results in the
literature. We
provide~\citep{FB-JC-SM:08cor,WR-RWB:08,MM-ME:10,ROS-JAF-RMM:07} and
the references therein as a starting point for the interested
reader. In cyberphysical systems, privacy at the physical layer
provides protection beyond the use of higher-level encryption-based
techniques.  Information-theoretic approaches to privacy at the
physical layer have been actively
pursued~\citep{DG-EE-HVP:10,AM-SAAF-JH-ALS:14}. Recently, these ideas
have also been utilized in the context of
control~\citep{TT-HS:15}. The paper~\citep{AM-SAAF-JH-ALS:14} also
surveys the more recent game-theoretic approach to the topic. In
computer science, the notion of differential privacy, first introduced
in~\citep{CD-FM-KN-AS:06,CD:06}, and the design of differentially
private mechanisms have been widely studied in the context of privacy
preservation of databases. The work~\citep{CD-AR:14} provides a recent
comprehensive treatment.  A well-known advantage of differential
privacy over other notions of privacy is its immunity to
post-processing and side information, which makes it particularly
well-suited for multi-agent scenarios where agents do not fully trust
each other and/or the communication channels are not fully secure.
While secure multi-party computation also deals with scenarios where
no trust exists among agents, the maximum number of agents that can
collude (without the privacy of others being breached) is bounded,
whereas using differential privacy provides immunity against arbitrary
collusions~\citep{PK-SO-PV:15, MP-PL:15}.  As a result, differential
privacy has been adopted by recent works in a number of areas
pertaining to networked systems, such as
control~\citep{ZH-SM-GD:12,ZH-YW-SM-GED:14,YW-ZH-SM-GED:14},
estimation~\citep{JLN-GJP:14}, and
optimization~\citep{SH-UT-GJP:14,ZH-SM-NV:15,EN-PT-JC:17-tcns}.  Of
relevance to our present work, the paper~\citep{ZH-SM-GD:12} studies
the average consensus problem with differentially privacy guarantees
and proposes an adjacency-based distributed algorithm with decaying
Laplace noise and mean-square convergence.  The algorithm preserves
the differential privacy of the agents' initial states but the
expected value of its convergence point depends on the network
topology and may not be the exact average, even in expectation.  By
contrast, the algorithm proposed in this work enjoys almost sure
convergence, asymptotic unbiasedness, and an explicit characterization
of its convergence rate.  Our results also allow individual agents to
independently choose their level of privacy.  The problem of
privacy-preserving average consensus has also been studied using other
notions of privacy.  The work~\citep{NEM-CNH:13} builds
on~\citep{MK-MST-BHK-HRR:07} to let agents have the option to add to
their first set of transmitted messages a zero-sum noise sequence with
finite random length, which in turn allows the coordination algorithm
to converge to the exact average of their initial states. As long as
an adversary cannot listen to the transmitted messages of an agent as
well as all its neighbors, the privacy of that agent is preserved, in
the sense that different initial conditions may produce the same
transmitted messages. This idea is further developed
in~\citep{YM-RMM:14,YM-RMM:15}, where agents add an infinitely-long
exponentially-decaying zero-sum sequence of Gaussian noise to their
transmitted messages. The algorithm has guaranteed mean-square
convergence to the average of the agents' initial states and preserves
the privacy of the nodes whose messages and those of their neighbors
are not listened to by the malicious nodes, in the sense that the
maximum-likelihood estimate of their initial states has nonzero
variance. Finally, \citep{XD-etal:15} considers the problem of privacy
preserving maximum consensus.

\emph{Statement of Contributions:} We study the average consensus
problem where a group of agents seek to compute and agree on the
average of their local variables while seeking to keep them
differentially private against an adversary with potential access to
all group communications.  This privacy requirement also applies to
the case where each agent wants to keep its initial state private
against the rest of the group (e.g., due to the possibility of
communication leakages).  The main contributions of this work are the
characterization and optimization of the fundamental trade-offs
between differential privacy and average consensus. Our first
contribution is the formulation and formal proof of a general
impossibility result. We show that as long as a coordination algorithm
is differentially private, it is impossible to guarantee the
convergence of agents' states to the average of their initial values,
even in distribution. This result automatically implies the same
impossibility result for stronger notions of convergence. Motivated by
it, our second contribution is the design of a linear Laplacian-based
consensus algorithm that achieves average consensus in expectation
---the most that one can expect. We prove the almost sure convergence
and differential privacy of our algorithm and characterize its
accuracy and convergence rate. Our final contribution is the
computation of the optimal values of the design parameters to achieve
the most accurate consensus possible. Letting the agents fix a (local)
desired value of the privacy requirement, we minimize the variance of
the algorithm convergence point as a function of the noise-to-state
gain and the amplitude and decay rate of the noise. We show that the
minimum variance is achieved by the one-shot perturbation of the
initial states by Laplace noise.  This result reveals the optimality
of one-shot perturbation for static average consensus, previously (but
implicitly) shown in the sense of information-theoretic entropy.
Various simulations illustrate our results.

\section{Preliminaries}\label{section:preliminaries}

This section introduces notations and basic concepts. We denote the set
of reals, positive reals, non-negative reals, positive integers, and
nonnegative integers by $\real$, $\realpos$, $\realnonneg$,
$\integerpos$, and $\integernonneg$, respectively.  We denote by~$\|
\cdot \|$ the Euclidean norm.  We let $(\real^n)^\integerpos$ denote
the space of vector-valued sequences in~$\real^n$. For $\lbrace x(k)
\rbrace_{k = 0}^\infty \in (\real^n)^\integerpos$, we use the
shorthand notation $\mathbf{x} = \lbrace x(k) \rbrace_{k = 0}^\infty$
and $\mathbf{x}_k = \lbrace x(j) \rbrace_{j = 0}^k$. $I_n \in \real^{n
  \times n}$ and $\ones_n \in \real^n$ denote the identity matrix and
the vector of ones, respectively. For $x \in \real^n$, $\Ave{x} =
\frac{1}{n} \ones_n^T x$ denotes the average of its components.  We
let $\Pi_n = \frac{1}{n} \ones_n \ones_n^T$. Note that $\Pi_n$ is
diagonalizable, has one eigenvalue equal to $1$ with eigenspace
\begin{align*}
  \real \ones_n \triangleq \lbrace a \ones_n \mid a \in \real
  \rbrace,
\end{align*}
and all other eigenvalues equal~$0$.  For a vector space $V \subset
\real^n$, we let $V^\perp$ denote the vector space orthogonal to $V$.
A matrix $A \in \real^{n \times n}$ is stable if all its eigenvalues
have magnitude strictly less than 1. A function $\gamma: [0, \infty)
\to [0, \infty)$ belongs to class $\Kc$ if it is continuous and
strictly increasing and $\gamma(0) = 0$. Similarly, a function
$\beta:[0, \infty) \times [0, \infty) \to [0, \infty)$ belongs to
class $\Kc \Lc$ if $\beta(\cdot, s)$ belongs to class $\Kc$ for any $s
\in [0, \infty)$ and $\beta(r, \cdot)$ is decreasing and $\lim_{s \to
  \infty} \beta(r, s) = 0$ for any $r \in [0, \infty)$. For $q \in
(0,1)$, the Euler function is given by $\varphi(q) = \prod_{k =
  1}^\infty (1 - q^k) > 0$. Note that
\begin{align*}
  \lim_{k \to \infty } \prod_{j = k}^\infty (1 - q^j) = \lim_{k \to
    \infty } \frac{\varphi(q)}{\prod_{j = 1}^{k - 1} (1 - q^j)} = 1.
\end{align*}
For a function $f:X \to Y$ and sets $A \subseteq X$ and $B \subseteq
Y$, we use $f(A) = \{f(x) \in Y | x \in A\}$ and $f^{-1}(B) = \{x \in
X | f(x) \in B\}$.  In general, $f(f^{-1}(B)) \subseteq B$. Finally,
for any topological space $X$, we denote by $\Bc(X)$ the set of Borel
subsets of~$X$.

\subsection{Graph Theory}

We present some useful notions on algebraic graph theory
following~\citep{FB-JC-SM:08cor}. Let $\Gc = (V,E,A)$ denote a
weighted undirected graph with vertex set $V$ of cardinality~$n$, edge
set $E \subset V \times V$, and symmetric adjacency matrix $A \in
\realnonneg^{n \times n}$. A path from $i$ to $j$ is a sequence of
vertices starting from $i$ and ending in $j$ such that any pair of
consecutive vertices is an edge. The set of neighbors $\Nc_i$ of $i$
is the set of nodes $j$ such that $(i,j) \in E$.  A graph is connected
if for each node there exists a path to any other node.  The weighted
degree matrix is the diagonal matrix $D \in \real^{n \times n}$ with
diagonal $A \ones_n$. The Laplacian is $L = D - A$ and has the
following properties:
\begin{itemize}
\item $L$ is symmetric and positive semi-definite;
\item $L \ones_n = 0$ and $\ones_n^T L = 0$, i.e., $0$ is an
  eigenvalue of $L$ corresponding to the eigenspace $\real \ones_n$;
\item $\Gc$ is connected if and only if $\text{rank}(L) = n-1$, so $0$
  is a simple eigenvalue of $L$;
\item All eigenvalues of $L$ belong to $[0, 2 d_\text{max}]$, where
  $d_\text{max}$ is the largest element of $D$.
\end{itemize}
For convenience, we define $L_\text{cpt} = I_n - \Pi_n$.

\subsection{Probability Theory}\label{subsection:probability}

Here we briefly review basic notions on probability
following~\citep{AP-SUP:02, RD:10}. Consider a probability space
$(\Omega, \Sigma, \pr)$.  If $E, F \in \Sigma$ are two events with $E
\subseteq F$, then $\pr \lbrace E \rbrace \le \pr \lbrace F \rbrace
$. For simplicity, we may sometimes denote events of the type $E_p =
\lbrace \omega \in \Omega \mid p(\omega) \rbrace$ by~$\lbrace p
\rbrace$, where $p$ is a logical statement on the elements of
$\Omega$.  Clearly, for two statements $p$ and $q$,
\begin{align}\label{eq:prob-implication}
  (p \Rightarrow q) \Rightarrow \left( \pr \lbrace p \rbrace \le \pr
    \lbrace q \rbrace \right).
\end{align}
A random variable is a measurable function $X:\Omega \to \real$. For
any $N \in \realpos$ and any random variable $X$ with finite expected
value $\mu$ and finite nonzero variance $\sigma^2$, Chebyshev's
inequality states that
\begin{align}\label{eq:Chebyshev}
  \pr \lbrace \lvert X - \mu \rvert \ge N \sigma \rbrace \le
  \frac{1}{N^2}.
\end{align}

For a random variable $X$, let $\E[X]$ and $F_X$ denote its
expectation and cumulative distribution function, respectively.  A
sequence of random variables $\lbrace X_k \rbrace_{k \in
  \integernonneg}$ converges to a random variable $X$
\begin{itemize}
\item almost surely (a.s.) if $\pr \lbrace \lim_{k \to \infty} X_k = X
  \rbrace = 1$;
\item in mean square if $\E[X_k^2], \E[X^2] < \infty$ for all $k \in
  \integernonneg$ and $\lim_{k \to \infty} \E[( X_k - X )^2] = 0$;
\item in probability if $\lim_{k \to \infty} \pr \lbrace \lvert X_k -
  X \rvert < \upsilon \rbrace = 1$ for any $\upsilon > 0$;
\item in distribution or weakly if $\lim_{k \to \infty} F_{X_k}(x) =
  F_X(x)$ for any $x \in \real$ at which $F_X$ is continuous.
\end{itemize}
Almost sure convergence and convergence in mean square imply
convergence in probability, which itself implies convergence in
distribution. Moreover, if $\pr \lbrace \lvert X_k \rvert \le \bar X
\rbrace = 1$ for all $k \in \integernonneg$ and some fixed random
variable $\bar X$ with $\E[\bar X^2] < \infty$, then convergence in
probability implies mean square convergence, and if $X$ is a constant,
then convergence in distribution implies convergence in probability.

A zero-mean random variable $X$ has Laplace distribution with scale $b
\in \realpos$, denoted $X \sim \Lap(b)$, if its pdf~is given by
$\Lc(x; b) \triangleq \frac{1}{2b} e^{-\frac{|x|}{b}}$ for any $x \in
\real$. It is easy to see that $\lvert X \rvert$ has exponential
distribution with rate $\lambda = \frac{1}{b}$.
  
\subsection{Input-to-State Stability of Discrete-Time
  Systems}\label{subsection:ISS}

This section briefly describes notions of robustness for discrete-time
systems following~\citep{ZPJ-YW:01}. Consider a discrete-time system
of the form
\begin{align}\label{eq:discrete-dynamics}
  x(k+1) = f(x(k),u(k)),
\end{align}
where $u: \integernonneg \to \real^m$ is a disturbance input, $x:
\integernonneg \to \real^n$ is the state, and $f: \real^n \times
\real^m \to \real^n$ is a vector field satisfying $f(0, 0) = 0$. The
system~\eqref{eq:discrete-dynamics} is globally input-to-state stable
(ISS) if there exists a class $\Kc \Lc$ function~$\beta$ and a class
$\Kc$ function~$\gamma$ such that, for any bounded input $u$,
 any initial condition $x_0 \in \real^n$, and all $k \in \integernonneg$,
\begin{align*}
  \|x(k)\| \le \beta(\|x_0\|, k) + \gamma(\|u\|_{\ell_\infty}),
\end{align*}
where $\|u\|_{\ell_\infty} = \sup \lbrace \|u(k)\| \mid k \in
\integernonneg \rbrace$. The system~\eqref{eq:discrete-dynamics} has a
$\Kc$-asymptotic gain if there exists a class $\Kc$ function
$\gamma_a$ such that, for any initial condition $x_0 \in \real^n$,
\begin{align*}
  \limsup_{k \to \infty} \|x(k)\| \le \gamma_a \big( \limsup_{k \to
    \infty} \|u(k)\| \big).
\end{align*}
If a system is ISS, then it has a $\Kc$-asymptotic gain. Furthermore,
any LTI system $x(k+1) = A x(k) + B u(k)$ is ISS if $A$ is stable.

\section{Problem statement}\label{section:problem-statement}

Consider a group of $n$ agents whose interaction topology is described
by an undirected connected graph~$\Gc$. The group objective is to
compute the average of the agents' initial states while preserving the
privacy of these values against potential adversaries eavesdropping
on all the network communications. Note that this privacy
  requirement is the same as the case where each agent wants to keep
  its initial state private against the rest of the group due to the
  possibility of communication leakages. We next generalize
the exposition in~\citep{ZH-SM-GD:12} to provide a formal presentation
of this problem.  The state of each agent $i \in \until{n}$ is
represented by $\theta_i \in \real $. The message that agent $i$
shares with its neighbors about its current state is denoted by~$x_i
\in \real$. For convenience, the aggregated network state and the
vector of transmitted messages are denoted by $\theta = (\theta_1 ,
\dots , \theta_n) \in \real^n$ and $x = (x_1 , \dots , x_n) \in
\real^n$, respectively.  Agents update their states in discrete time
according to some rule,
\begin{align}\label{eq:dynamics-nonlin}
  \theta(k+1) = f (\theta(k), x(k) ), \quad k \in \integernonneg ,
\end{align}
with initial states $\theta(0) = \theta_0$, where the state-transition
function $f:\real^n \times \real^n \to \real^n$ is such that its $i$th
element depends only on $\theta_i$ and $\{x_j\}_{j \in \Nc_i \cup
  \lbrace i \rbrace }$. The messages are calculated as
\begin{align}\label{eq:message-nonlin}
  x(k) = h(\theta(k), \eta(k)), \quad k \in \integernonneg,
\end{align}
where $h: \real^n \times \real^n \to \real^n$ is such that its $i$th
element depends only on $\theta_i$ and $\eta_i$. For simplicity,
  we assume that $f$ and $h$ are continuous. $\eta(k) \in \real^n$ is
a vector random variable, with $\eta_i(k)$ being the noise generated
by agent $i$ at time $k$ from an arbitrary distribution. Consequently,
$\boldsymbol{\theta}$ and $\mathbf{x}$ are sequences of vector random
variables on the total sample space $\Omega = (\real^n)^\integerpos$
whose elements are noise sequences $\boldsymbol{\eta}$. Although one
could choose $h$ to only depend on $\theta$, corrupting the messages
by noise is necessary to preserve privacy.  Given an initial state
$\theta_0$, $\mathbf{x}$ is uniquely determined by $\boldsymbol{\eta}$
according
to~\eqref{eq:dynamics-nonlin}-\eqref{eq:message-nonlin}. Therefore,
the function $X_{\theta_0}: (\real^n)^\integerpos \to
(\real^n)^\integerpos$ such that
\begin{align*}
  X_{\theta_0} (\boldsymbol{\eta}) = \mathbf{x}
\end{align*}
is well defined.

\begin{definition}\longthmtitle{Differential
    Privacy}\label{def:privacy}
    Given $\delta \in \realpos$, the initial network states
    $\theta_0^{(1)}$ and $\theta_0^{(2)}$ are $\delta$-adjacent if,
    for some $i_0 \in \lbrace 1 , \dots , n \rbrace $,
  \begin{align}\label{eq:adjacency}
    \lvert \theta_{0,i}^{(2)} - \theta_{0,i}^{(1)} \rvert \le
    \begin{cases}
      \delta & \text{if $i = i_0$,}
      \\
      0 & \text{if $i \neq i_0$,}
    \end{cases}
  \end{align}
  for $ i \in \until{n}$.  Given $\delta, \epsilon \in \realnonneg$,
  the dynamics~\eqref{eq:dynamics-nonlin}-\eqref{eq:message-nonlin} is
  $\epsilon$-differentially private if, for any pair $\theta_0^{(1)}$
  and $\theta_0^{(2)}$ of $\delta$-adjacent initial states and any set
  $\Ocal \in \Bc\big((\real^n)^\integerpos\big)$,
  \begin{align*}
    \pr \lbrace \boldsymbol{\eta} \in \Omega \! \mid \! X_{\theta_0^{(1)}}
    (\boldsymbol{\eta}) \in \Ocal \rbrace \le e^{\epsilon} \pr \lbrace
    \boldsymbol{\eta} \in \Omega \mid X_{\theta_0^{(2)}}
    (\boldsymbol{\eta}) \in \Ocal \rbrace .
  \end{align*}
\end{definition}

A final aspect to consider is that, because of the presence of noise,
the agents' states under~\eqref{eq:dynamics-nonlin} might not converge
exactly to their initial average $\Ave{\theta_0}$, but to a
neighborhood of it. This is captured by the notion of accuracy.

\begin{definition}\longthmtitle{Accuracy}\label{def:accuracy}
  For $p \in [0,1]$ and $r \in \realnonneg$, the
  dynamics~\eqref{eq:dynamics-nonlin}-\eqref{eq:message-nonlin} is
  $(p,r)$-accurate if, from any initial state $\theta_0$, the network
  state $\theta(k)$ converges to $\theta_\infty \in \real^n$ as $k \to
  \infty$, with $ \E[\theta_\infty] = \Ave{\theta_0} \ones_n$ and $
  \pr \lbrace \| \theta_\infty - \Ave{\theta_0} \ones_n \| \le r
  \rbrace \ge 1 - p$.
\end{definition}

In Definition~\ref{def:accuracy}, the type of convergence of
$\theta(k)$ to $\theta_\infty$ can be any of the four classes
described in Section~\ref{subsection:probability}. Furthermore, for
each notion of convergence, $(0,0)$-accuracy is equivalent to the
convergence of $\theta(k)$ to $\Ave{\theta_0} \ones_n$.  We are
finally ready to formally state our problem.

\begin{problem}\longthmtitle{Differentially Private Average
    Consensus}\label{problem:main}
  Design the dynamics~\eqref{eq:dynamics-nonlin}, the inter-agent
  messages~\eqref{eq:message-nonlin}, and the distribution of noise
  sequences $\boldsymbol{\eta}$ such that asymptotic average consensus
  is achieved with $(p,r)$-accuracy while guaranteeing
  $\epsilon$-differential privacy for (finite) $\epsilon$, $r$, and $p
  \in \realnonneg$ as small as possible. \oprocend
\end{problem}

\section{Obstructions to Exact Differentially Private Average
  Consensus}\label{section:impossibility}

In this section we establish the impossibility of solving
Problem~\ref{problem:main} with $(0,0)$-accuracy, even if considering
the weakest notion of convergence.

\begin{proposition}[Impossibility Result]\label{prop:bound}
  Consider a group of agents executing a distributed algorithm of the
  form~\eqref{eq:dynamics-nonlin} with messages generated according
  to~\eqref{eq:message-nonlin}.  Then, for any $\delta, \epsilon > 0$,
  agents cannot simultaneously converge to the average of their
  initial states in distribution and preserve $\epsilon$-differential
  privacy of their initial states.
\end{proposition}
\begin{pf}
  We reason by contradiction. Assume there exists an algorithm that
  achieves convergence in distribution to the exact average of the
  network initial state and preserves $\epsilon$-differential privacy
  of it. Since the algorithm must preserve the privacy of \textit{any}
  pair of $\delta$-adjacent initial conditions, consider a specific
  pair satisfying
  \begin{align*}
    \theta_{0, i_0}^{(2)} = \theta_{0, i_0}^{(1)} + \delta,
  \end{align*}
  for some $i_0 \in \until{n}$ and $\theta_{0, i}^{(2)} = \theta_{0,
    i}^{(1)}$ for all $i \neq i_0$. Since $\Ave{\theta_0}$ is fixed
  (i.e., deterministic), the convergence of $\theta_i(k)$, $i
  \in \until{n}$ to $\Ave{\theta_0}$ is also in probability. Thus, for
  any $i \in \until{n}$ and any $\upsilon > 0$, we have $ \lim_{k \to
    \infty} \pr \lbrace \lvert \theta_i^{(\ell)} (k) -
  \Ave{\theta_0^{(\ell)}} \rvert < \upsilon \rbrace = 1$, for $ \ell =
  1,2$.  Therefore, for any $\upsilon' > 0$, there exists $k \in
  \integernonneg$ such that for all $i \in \until{n}$,
  \begin{align}\label{eq:prob-of-neighborhood}
    \pr \lbrace \lvert \theta_i^{(\ell)} (k) - \Ave{\theta_0^{(\ell)}}
    \rvert < \upsilon \rbrace > 1 - \upsilon', \quad \ell = 1,2
    .
  \end{align}
  Now, considering~\eqref{eq:dynamics-nonlin}-\eqref{eq:message-nonlin}, it
  is clear that, for any fixed initial state $\theta_0$ and any $k \in
  \integernonneg$, $\mathbf{x}_k$ is uniquely determined by
  $\boldsymbol{\eta}_k$ and $\boldsymbol{\theta}_k$ is uniquely
  determined by $\mathbf{x}_k$. Therefore, the functions $X_{k,
    \theta_0}, \Theta_{k, \theta_0} : \real^{n(k + 1)} \to \real^{n(k
    + 1)}$ such that
  \begin{align}\label{eq:mappings}
    X_{k, \theta_0} (\boldsymbol{\eta}_k) = \mathbf{x}_k, \quad
    \Theta_{k, \theta_0} (\mathbf{x}_k) = \boldsymbol{\theta}_k ,
  \end{align}
  are well defined and continuous (due to continuity of $f$ and
    $g$). Next, for $\ell = 1, 2$, define $R_k^{(\ell)} = X_{k,
      \theta_0^{(\ell)}}^{-1}\big(\Theta_{k,
      \theta_0^{(\ell)}}^{-1}\big(\Nc_k^{(\ell)}\big)\big)$, where
    $\Nc_k^{(\ell)} \triangleq \real^{n k} \times
    \big(\Ic^{(\ell)}\big)^n$ and $\Ic^{(\ell)} \subset \real$ is the
    $\upsilon$-neighborhood of $\Ave{\theta_0^{(\ell)}}$.
    By~\eqref{eq:prob-of-neighborhood}, we have
    \begin{align}\label{eq:prob-of-R_k^1}
      \pr (R_k^{(\ell)}) > 1 - \upsilon', \quad \ell = 1, 2.
    \end{align}
    Note that $R_k^{(1)}$ is open as it is the continuous pre-image of
    an open set, so $\Ocal_k \triangleq X_{k, \theta_0^{(1)}}
    \big(R_k^{(1)}\big)$ is Borel.  To reach a contradiction, we
    define $R_k^{'(2)} = X_{k, \theta_0^{(2)}}^{-1} (\Ocal_k)$ and
    show that $\pr (R_k^{'(2)})$ can be made arbitrarily small by
    showing that $R_k^{'(2)} \cap R_k^{(2)} = \emptyset$. To do this,
    note that by the definitions of $R_k^{'(2)}$, $\Ocal_k$ and
    $R_k^{(1)}$, we have
    \begin{align}\label{eq:Nk1}
      \Theta_{k, \theta_0^{(1)}}\big(X_{k,
        \theta_0^{(2)}}\big(R_k^{'(2)}\big)\big) \subseteq
      \Nc_k^{(1)}.
    \end{align}
    Recall that in~\eqref{eq:dynamics-nonlin}, $f$ is such that the
    next state of each agent only depends on its current state and the
    messages it receives. Hence, since for all $i \neq i_0$,
    $\theta_{0, i}^{(2)} = \theta_{0, i}^{(1)}$, we have
    from~\eqref{eq:Nk1} that
    \begin{align*}
      \Theta_{k, \theta_0^{(2)}}\big(X_{k,
        \theta_0^{(2)}}\big(R_k^{'(2)}\big)\big) \subseteq \overline
      \Nc_k^{(1)},
    \end{align*}
    where $\overline \Nc_k^{(1)} \triangleq \real^{n k} \times
    (\Ic^{(1)})^{i_0 - 1} \times \real \times (\Ic^{(1)})^{n - i_0}$
    is the same as $\Nc_k^{(1)}$ except that the requirement on
    $\theta_{i_0}(k)$ (to be close to $\Ave{\theta_0^{(1)}}$) is relaxed.
       Now, since $\Theta_{k, \theta_0^{(2)}}\big(X_{k,
      \theta_0^{(2)}}\big(R_k^{(2)}\big)\big) \subseteq \Nc_k^{(2)}$
    and, by choosing $\upsilon < \frac{\delta}{2n}$, we get $\overline
    \Nc_k^{(1)} \cap \Nc_k^{(2)} = \emptyset$, we conclude that
    $\Theta_{k, \theta_0^{(2)}}\big(X_{k,
      \theta_0^{(2)}}\big(R_k^{(2)}\big)\big) \cap \Theta_{k,
      \theta_0^{(2)}}\big(X_{k,
      \theta_0^{(2)}}\big(R_k^{'(2)}\big)\big) = \emptyset$, which
    implies $R_k^{(2)} \cap R_k^{'(2)} = \emptyset$,
    so we get
    \begin{align}\label{eq:prob-of-R_k^2}
      \pr (R_k^{(2)}) < \upsilon',
    \end{align}
    as desired. Now, let $\Ocal = \Ocal_k \times
      (\real^n)^\integerpos \in
      \Bc\big((\real^n)^\integerpos\big)$. For any initial condition
    $\theta_0$,
    \begin{align*}
      \pr \{\boldsymbol{\eta} | X_{\theta_0}(\boldsymbol{\eta}) \in
      \Ocal \} = \pr \{\boldsymbol{\eta}_k | X_{k,
        \theta_0}(\boldsymbol{\eta_k}) \in \Ocal_k\}.
    \end{align*}
    Hence, since the algorithm is $\epsilon$-differentially private,
    \begin{align*}
      \pr (R_k^{(1)}) = & \pr \{\boldsymbol{\eta}_k | X_{k,
        \theta_0^{(1)}}(\boldsymbol{\eta_k}) \in \Ocal_k\}
      \\
      &\le e^\epsilon \pr \{\boldsymbol{\eta}_k | X_{k,
        \theta_0^{(2)}}(\boldsymbol{\eta_k}) \in \Ocal_k\} =
      e^\epsilon \pr (R_k^{'(2)}).
    \end{align*}
  Thus, using~\eqref{eq:prob-of-R_k^1} and~\eqref{eq:prob-of-R_k^2},
  we have for all $\upsilon' > 0$,
  \begin{align*}
    1 - \upsilon' < e^\epsilon \upsilon' \Rightarrow \frac{1}{1
      + e^\epsilon} < \upsilon' ,
  \end{align*}
  which is clearly a contradiction because $\epsilon$ is a finite
  number, completing the proof.  \qed
\end{pf}

Since convergence in distribution is the weakest notion of
convergence, Proposition~\ref{prop:bound} implies that a differentially
private algorithm cannot guarantee any type of convergence to the
exact average. Therefore, in our forthcoming discussion, we relax the
exact convergence requirement and allow for convergence to a random
variable that is at least unbiased (i.e., centered at the true
  average).

\section{Differentially Private Average Consensus
  Algorithm}\label{section:dynamics-design-and-analysis}

Here, we develop a solution to Problem~\ref{problem:main}. Consider
the following linear distributed dynamics,
\begin{align}\label{eq:dynamics}
  \theta(k+1) = \theta(k) - h L x(k) + S \eta(k) ,
\end{align}
for $ k \in \integernonneg$, where $h < (d_\text{max})^{-1}$ is the
step size, $S$ is a diagonal matrix with diagonal $(s_1, \dots, s_n)$
and $s_i \in (0,2)$ for each $i \in \until{n}$, and the
messages are generated as
\begin{align}\label{eq:message}
  x(k) = \theta (k) + \eta(k),
\end{align}
where the $i$th component of the noise vector $\eta(k)$ has the
  Laplace distribution $\eta_i(k) \sim \Lap(b_i(k))$ at any time $k
  \in \integernonneg$ with
  \begin{align}\label{eq:noise}
    b_i(k) = c_i q_i^k, \quad c_i \in \realpos, \; q_i \in (\lvert s_i
    - 1 \rvert, 1) .
  \end{align}
Note that~\eqref{eq:dynamics} is a special case
of~\eqref{eq:dynamics-nonlin} (since $\eta(k) = x(k) - \theta(k)$)
and~\eqref{eq:message} a special case
of~\eqref{eq:message-nonlin}. Also note that without the term $S
\eta(k)$, the average of the agents' initial states would be preserved
throughout the evolution.

  \begin{remark}\longthmtitle{Comparison with the Literature}
    {\rm The  proposed algorithm~\eqref{eq:dynamics}-\eqref{eq:noise}
      has similarities and differences with the algorithm proposed
      in~\citep{ZH-SM-GD:12} which can be expressed  (with a slight
      change of notation in using $s_i$ instead of $\sigma_i$) as
      \begin{align*}
        \theta(k+1) &= (I_n - S) \theta(k) + S D^{-1} A x(k)
        \\
        &= \big[I_n - S D^{-1} L \big] \theta(k) + \big[S - S D^{-1} L
        \big] \eta(k).
      \end{align*}
      If each agent selects $s_i = d_i h < 1$, then we
      recover~\eqref{eq:dynamics}-\eqref{eq:noise}. As we show later,
      this particular choice results in an unbiased convergence point,
      while in general the expected value of the convergence point of
      the algorithm in~\citep{ZH-SM-GD:12} depends on the graph
      structure and may not be the true average.  Furthermore, this
      algorithm is only shown to exhibit mean square convergence of
      $\theta(k)$ for $s_i \in (0, 1)$, while here we provide an
      explicit expression for the convergence point and establish
      convergence in the stronger a.s. sense for larger range of $s_i
      \in (0, 2)$.  As we show later, the inclusion of $s_i = 1$ is
      critical, as it leads to identifying the optimal algorithm
      performance.  On a different note, the algorithms
      in~\citep{YW-ZH-SM-GED:14} and~\citep{YM-RMM:14,YM-RMM:15} add a
      noise sequence to the messages which is correlated over time --
      the latter using a different notion of
      privacy. \cite{YW-ZH-SM-GED:14} generate a single noise at time
      $k = 0$ and add a scaled version of it to the messages at every
      time $k \ge 1$, leading to an effectively ``one-shot''-type of
      perturbation. We show in
      Section~\ref{section:optimal-noise-selection} that the one-shot
      approach is optimal for static average consensus while
      sequential perturbation is necessary for dynamic scenarios.
      \oprocend}
\end{remark}

\subsection{Convergence Analysis}

This section analyzes the asymptotic correctness of the
algorithm~\eqref{eq:dynamics}-\eqref{eq:noise} and characterizes its rate of
convergence. We start by establishing convergence.

\begin{proposition}\longthmtitle{Asymptotic
    Convergence}\label{prop:convergence} 
  Consider a network of $n$ agents executing the distributed
  dynamics~\eqref{eq:dynamics}-\eqref{eq:noise}. Define the random
  variable $\theta_\infty$ as
  \begin{align}\label{eq:theta-infty}
    \theta_\infty \triangleq \Ave{\theta_0} + \sum_{i = 1}^n \frac{s_i}{n}
    \sum_{j = 0}^\infty \eta_i(j) .
  \end{align}
  Then, $\theta_\infty$ is well-defined a.s., and the states of all
  agents converge to $\theta_\infty$ almost surely.
\end{proposition}
\begin{pf}
  Note that $s_i \in (0,2)$ ensures that $(\lvert s_i - 1 \rvert, 1)$
  is not empty.  Substituting~\eqref{eq:message}
  into~\eqref{eq:dynamics}, the system dynamics is
  \begin{align}\label{eq:dynamics-without-x}
    \theta(k+1) = A \theta(k) + B \eta(k) ,
  \end{align}
  with $A = I_n - h L$ and $B=S - h L$. For any $\theta \in \real^n$,
  let
  \begin{align}\label{eq:theta-tilde}
    \tilde \theta = \theta - \Ave{\theta} \ones_n = L_\text{cpt}
    \theta \in (\real \ones_n)^\perp.
  \end{align}
  Multiplying both sides of~\eqref{eq:dynamics-without-x} by
  $L_\text{cpt}$ on the left and using the fact that $L_\text{cpt}$
  and $L$ commute, the dynamics of $\tilde \theta$~can be expressed as
  \begin{align}\label{eq:theta-tilde-dynamics}
    \tilde \theta(k+1) = (I_n - h L) \tilde \theta(k) + L_\text{cpt} (S
    - h L) \eta(k) .
  \end{align}
  Notice that $(\real \ones_n)^\perp$ is forward invariant
  under~\eqref{eq:theta-tilde-dynamics}. Therefore, considering
  $(\real \ones_n)^\perp$ as the state space for~\eqref{eq:theta-tilde-dynamics}
  and noting that $I_n - h L$ is stable on it, we deduce
  that~\eqref{eq:theta-tilde-dynamics} is ISS. Consequently, this
  dynamics has a $\Kc$-asymptotic gain
  (c.f. Section~\ref{subsection:ISS}), i.e., there exists $\gamma_a
  \in \Kc$ such that
  \begin{align*}
    \limsup_{k \to \infty} \|\tilde \theta(k)\| \le \gamma_a \big
    (\limsup_{k \to \infty} \|\eta(k)\| \big).
  \end{align*}
  Therefore, $\lim_{k \to \infty} \tilde \theta(k) \neq 0$ implies
  $\lim_{k \to \infty} \|\eta(k)\| \neq 0$. By definition, the latter
  means that there is $\upsilon > 0$ such that for all $K \in
  \integerpos$ there exists $k \ge K$ with $ \|\eta(k)\| > \upsilon$.
  In other words, there exists a subsequence $\left \lbrace
    \eta(k_\ell) \right \rbrace_{\ell \in \integerpos}$ such that
  $\|\eta(k_\ell)\| > \upsilon$ for all $\ell \in \integerpos$. This,
  in turn, implies that for all $\ell \in \integerpos$, $
  \|\eta(k_\ell)\|_\infty > \upsilon / \sqrt{n}$, i.e.,
  \begin{align*}
    \exists i_\ell \in \until{n} \text{ with } \lvert \eta_{i_\ell}
    (k_\ell) \rvert > \frac{\upsilon}{\sqrt{n}} .
  \end{align*}
  According to~\eqref{eq:prob-implication}, this chain of implications
  gives
  \begin{align*}
    \pr \lbrace \lim_{k \to \infty} \tilde \theta(k) \neq 0 \rbrace
    &\le \pr \lbrace \forall \ell \in \integerpos, \; \exists i_\ell
    \text{ s.t. } \lvert \eta_{i_\ell}(k_\ell) \rvert >
    \frac{\upsilon}{\sqrt{n}} \rbrace
    \\
    &= \prod_{\ell = 1}^\infty e^{-\frac{\upsilon}{\sqrt{n}
        b_{i_\ell}(k_\ell)}} = 0.
  \end{align*}
  The last equality holds because $\lim_{\ell \to \infty}
  b_{i_\ell}(k_\ell) = \lim_{\ell \to \infty} c_{i_\ell}
  q_{i_\ell}^{k_\ell} = 0$. Therefore, we conclude
  \begin{align}\label{eq:as-1}
    \pr \lbrace \lim_{k \to \infty} \tilde \theta(k) = 0 \rbrace = 1.
  \end{align}
  From~\eqref{eq:theta-tilde}, we see that a.s. convergence of
  $\theta$ requires a.s. convergence of $\Ave{\theta}$ as well. Left
  multiplying~\eqref{eq:dynamics} by $\ones_n^T$, we obtain for all $k
  \in \integernonneg$,
  \begin{align*}
    \nonumber \frac{1}{n} \ones_n^T \theta(k+1) &= \frac{1}{n}
    \ones_n^T \theta(k) + \frac{1}{n} \ones_n^T S \eta(k)
    \\
    \nonumber &= \frac{1}{n} \ones_n^T \theta_0 + \frac{1}{n}
    \sum_{j=0}^k \sum_{i = 1}^n s_i \eta_i(j) ,
  \end{align*}
  which in turn implies
  \begin{align}\label{eq:dynamics-ave}
    \Ave{\theta(k)} &= \Ave{\theta_0} + \sum_{i = 1}^n \frac{s_i}{n}
    \sum_{j = 0}^{k - 1} \eta_i(j).
  \end{align}
  We next prove that $\Ave{\theta(k)}$ converges almost surely
  to~$\theta_\infty$. For the latter to be well-defined over~$\Omega$,
  we simply set $\theta_\infty \triangleq \Ave{\theta_0}$ when the
  series does not converge.  Clearly, for any $\boldsymbol{\eta} \in
  \Omega$ such that $\sum_{j = 0}^\infty \eta_i(j)$ converges for all
  $i \in \until{n}$, we have $\lim_{k \to \infty} \Ave{\theta(k)} =
  \theta_\infty$. Hence, using~\eqref{eq:prob-implication},
  \begin{align*}
    \pr \lbrace \lim_{k \to \infty} \Ave{\theta(k)} = \theta_\infty
    \rbrace \ge \prod_{i = 1}^n \pr \Big \{ \sum_{j = 0}^\infty
    \eta_i(j) \ \text{converges} \Big \} .
  \end{align*}
  Note that, for each $i \in \until{n}$ and any $\ell \in \integerpos$,
  if $ |\eta_i(j)| \le \frac{1}{j^2}$ for all $j \ge \ell$, then $
  \sum_{j = 0}^\infty \eta_i(j)$ converges.  Hence,
  using~\eqref{eq:prob-implication} and the definition of Laplace
  distribution, we get for all $\ell \in \integerpos$,
  \begin{align*}
    \pr \lbrace \lim_{k \to \infty} \Ave{\theta(k)} = \theta_\infty
    \rbrace & \ge \prod_{i = 1}^n \prod_{j = \ell}^\infty \pr \Big\{
    |\eta_i(j)| \le \frac{1}{j^2} \Big \}
    \\
    &= \prod_{i = 1}^n \prod_{j = \ell}^\infty \Big(1 - e^{-\frac{1}{c_i
        q_i^j j^2}} \Big).
  \end{align*}
  For each $i \in \until{n}$, because $0 < q_i < 1$, there exists $
  \beta_i$ such that $\frac{1}{c_i q_i^j j^2} \ge \beta_i j$ for $j
  \ge 1$. Therefore, using the Euler function $\varphi$,
  \begin{align*}
    \pr \lbrace \lim_{k \to \infty} \Ave{\theta(k)} & = \theta_\infty
    \rbrace \ge \prod_{i = 1}^n \frac{\varphi(e^{-\beta_i})}{\prod_{j
        = 1}^{\ell - 1} (1 - e^{-\beta_i j})},
 \end{align*}
 for all $\ell \in \integerpos$, and hence,
 \begin{align*}
   \pr \lbrace \lim_{k \to \infty} \Ave{\theta(k)} = &\ \theta_\infty
   \rbrace \!\ge \lim_{\ell \to \infty} \prod_{i = 1}^n
   \frac{\varphi(e^{-\beta_i})}{\prod\limits_{j = 1}^{\ell - 1} (1 -
     e^{-\beta_i j})} = 1.
 \end{align*}
 This, together with~\eqref{eq:theta-tilde} and~\eqref{eq:as-1},
 implies that $ \pr \lbrace \lim_{k \to \infty} \theta(k) = \theta_\infty
 \ones_n \rbrace = 1$, which completes the proof. \qed
\end{pf}

 \begin{remark}[Mean-Square Convergence]\label{rem:ms}
  {\rm From~\eqref{eq:dynamics-without-x} and the fact that $\|A\| =
    1$, we have
      \begin{align*}
        \|\theta(k)\| &\le \|\theta_0\| + \|B\| \sum_{j = 0}^{k-1}
        \|\eta(j)\|
        \\
        &\le \|\theta_0\| + \|B\| \sum_{j = 0}^\infty \|\eta(j)\|
        \triangleq Z ,
      \end{align*}
      for all $k \in \integernonneg$. It is straightforward to show
      $\E[Z^2] < \infty$, so, using
      Proposition~\ref{prop:convergence}, $\theta(k)$ also converges
      to $\theta_\infty \ones_n$ in mean square.}  \oprocend
\end{remark}

Our next aim is to characterize the convergence rate of the
distributed dynamics~\eqref{eq:dynamics}-\eqref{eq:noise}. Given the
result in Proposition~\ref{prop:convergence}, we define the exponential
mean-square convergence rate of the
dynamics~\eqref{eq:dynamics}-\eqref{eq:noise} as
\begin{align*}
  \hspace{-0.5ex}\mu = \!\lim_{k \to \infty} \!\left(\!\sup_{
      \theta(0) \in \real^n}
    \frac{\E\big[(\theta(k) - \theta_\infty \ones_n)^T (\theta(k) -
      \theta_\infty \ones_n) \big]}{\E\big[(\theta(0) - \theta_\infty
      \ones_n)^T (\theta(0) - \theta_\infty \ones_n)\big]}
  \right)^\frac{1}{2k} \hspace*{-1ex}.
\end{align*}
In the absence of noise ($\boldsymbol{\eta} = 0$), this definition
coincides with the conventional exponential convergence rate of
autonomous linear systems, see e.g.,~\cite{FB-JC-SM:08cor}.

\begin{proposition}\longthmtitle{Convergence
    Rate}\label{prop:conv-rate}
  Under the hypotheses of Proposition~\ref{prop:convergence}, the
  exponential mean-square convergence rate of the distributed
  dynamics~\eqref{eq:dynamics}-\eqref{eq:noise}~is
  \begin{align}\label{eq:mu}
    \mu = \max\{\overline q, \overline \lambda\} \in (0, 1),
  \end{align}
  where $\displaystyle \overline q = \max_{1 \le i \le n} q_i$ and
  $\overline \lambda < 1$ is the spectral radius of $I_n - h L -
  \Pi_n$.
\end{proposition}
\begin{pf}
  For convenience, we let $\hat \theta(k) = \theta(k) - \theta_\infty
  \ones_n$ denote the convergence error at $k \in \integernonneg$ and
  $\hat \theta_0 = \hat \theta(0)$.  Our first goal is to obtain an
  expression for $\E\big[\hat \theta(k)^T \hat \theta(k)\big]$.
  From~\eqref{eq:theta-infty} and the proof of
  Proposition~\ref{prop:convergence}, we have
  \begin{align*}
    \theta_\infty = \frac{1}{n} \ones_n^T \theta_0 + \frac{1}{n}
    \ones_n^T S \sum_{j = 0}^\infty \eta(j),
  \end{align*}
  almost surely.  Then, from~\eqref{eq:dynamics-without-x}, we have
  almost surely for all $k \in \integernonneg$,
  \begin{align*}
    \hat \theta(k) = A^k \theta_0 &+ \sum_{j = 0}^{k - 1} A^{k - 1 -
      j} B \eta(j) - \Pi_n \theta_0 - \Pi_n B \sum_{j = 0}^\infty
    \eta(j),
  \end{align*}
  where we have used the fact that $\Pi_n S = \Pi_n B$. Next, note
  that for all $k \in \integerpos$,
  \begin{align}\label{eq:A-Pi_n-k}
    (A - \Pi_n)^k &= \sum_{j = 0}^k {k \choose j} (-\Pi_n)^{k - j} A^j
    \\
    &= A^k + \sum_{j = 0}^{k - 1} {k \choose j} (-1)^{k - j} \Pi_n =
    A^k - \Pi_n, \notag
  \end{align}
  where we have used the facts that $\Pi_n$ is idempotent and $\Pi_n
  A^j = \Pi_n$ for any $j \in \integernonneg$. Let $\Ac = A -
  \Pi_n$. Notice that $\Ac$ has spectral radius $\overline \lambda <
  1$ and the same eigenvectors as $L$. Then, using~\eqref{eq:A-Pi_n-k}
  twice, we have almost surely for all $k \in \integerpos$,
  \begin{align*}
    \hat \theta(k) = \Ac^k \theta_0 &+ \sum_{j = 0}^{k - 2} \Ac^{k - 1
      - j} B \eta(j)
    \\
    &+ L_\text{cpt} B \eta(k - 1) - \sum_{j = k}^\infty \Pi_n S
    \eta(j).
  \end{align*}
  By the independence of $\{\eta(j)\}_{j = 0}^\infty$ over time, we
  have
  \begin{align}\label{eq:e-th2}
    \notag \E\big[\hat \theta(k)^T \hat \theta(k)\big] &= \theta_0^T
    \Ac^{2k} \theta_0
    \\
    \notag &\quad + \sum_{j = 0}^{k - 2} \E[\eta(j)^T B \Ac^{2k - 2 -
      2j} B \eta(j)]
    \\
    \notag &\quad + \E[\eta(k - 1)^T B L_\text{cpt}^2 B \eta(k - 1)]
    \\
    &\quad + \sum_{j = k}^\infty \E[\eta(j)^T S \Pi_n^2 S \eta(j)],
  \end{align}
  for all $k \in \integerpos$. Next, we upper bound the exponential
  mean-square convergence rate~$\mu$. Let $\displaystyle \overline c =
  \max_{1 \le i \le n} c_i$ and note that for any $N \in \real^{n
    \times n}$ and any $j \in \integernonneg$,
  \begin{align*}
    \E[\eta(j)^T N^T N \eta(j)] &= \sum_{i = 1}^n 2 b_i^2(j) (N^T
    N)_{ii} \\
    &\le 2 \overline c^2 \overline q^{2j} \tr(N^T N) = 2 \overline c^2
    \overline q^{2j} \|N\|_F^2,
  \end{align*}
  where $\|\cdot\|_F$ denotes the Frobenius norm. Therefore,
  \begin{align*}
    \E\big[\hat \theta(k)^T &\hat \theta(k)\big] \le \theta_0^T
    \Ac^{2k} \theta_0 + 2 \overline c^2 \sum_{j = 0}^{k - 2} \overline
    q^{2j} \|\Ac^{k - 1 - j} B\|_F^2
    \\
    &+ 2 \overline c^2 \overline q^{2(k - 1)} \|L_\text{cpt} B\|_F^2 +
    2 \overline c^2 \sum_{j = k}^\infty \overline q^{2j} \|\Pi_n
    S\|_F^2.
  \end{align*}
  Since the Frobenius norm is submultiplicative, $\|N\|_F^2 \le n
  \|N\|^2$ for any matrix $N$, and $\|\Ac\| = \overline \lambda$, we
  have
  \begin{align*}
    \E\big[\hat \theta(k)^T &\hat \theta(k)\big] \le \theta_0^T
    \Ac^{2k} \theta_0 + C_1 \sum_{j = 0}^{k - 2} \overline q^{2j}
    \overline \lambda^{2k - 4 - 2j} + C_2 \overline q^{2k} ,
  \end{align*}
  where $C_1 = 2 n \overline c^2 \|B\|_F^2 \overline \lambda^2$ and
  $C_2 = 2 \overline c^2(\|L_\text{cpt} B\|_F^2 / \overline q^2 +
  \|\Pi_n S\|_F^2 / (1 - \overline q^2))$ are constants. Note that for
  any $0 \le j \le k - 2$, we have $\overline q^{2j} \overline
  \lambda^{2k - 4 - 2j} \le \max\{\overline q, \overline \lambda\}^{2k
    - 4}$. Therefore, using the fact that the supremum of a sum is
  less than or equal to the sum of suprema, we have
  \begin{align*}
    &\sup_{\theta_0 \in \real^n} \frac{\E\big[\hat \theta(k)^T \hat
      \theta(k)\big]}{\E\big[\hat \theta_0^T \hat \theta_0\big]} \le
    \sup_{\theta_0 \in \real^n} \frac{\theta_0^T \Ac^{2k}
      \theta_0}{\E\big[\hat \theta_0^T \hat \theta_0\big]}
    \\
    &\qquad \qquad \qquad + \sup_{\theta_0 \in \real^n}\frac{C_3 (k -
      1) \max\{\overline q, \overline \lambda\}^{2k} + C_2 \overline
      q^{2k}}{\E\big[\hat \theta_0^T \hat \theta_0\big]},
  \end{align*}
  where $C_3 = C_1 \max\{\overline q, \overline \lambda\}^{-4}$. Let
  $\tilde \theta_0 = L_\text{cpt} \theta_0$ be the initial
  disagreement vector. It is straightforward to verify that
  $\theta_0^T \Ac^{2k} \theta_0 = \tilde \theta_0^T \Ac^{2k} \tilde
  \theta_0$ and
  \begin{align*}
    \E\big[\hat \theta_0^T \hat \theta_0\big] = \tilde \theta_0^T
    \tilde \theta_0 + \frac{1}{n} \sum_{i = 1}^n \frac{2 c_i^2
      s_i^2}{1 - q_i^2} \triangleq \tilde \theta_0^T \tilde \theta_0 +
    C_4.
  \end{align*}
  Therefore,
  \begin{align*}
    \sup_{\theta_0 \in \real^n} \frac{\E\big[\hat \theta(k)^T \hat
      \theta(k)\big]}{\E\big[\hat \theta_0^T \hat \theta_0\big]} &\le
    \sup_{\tilde \theta_0 \in (\real \ones_n)^\perp} \frac{\tilde \theta_0^T
      \Ac^{2k} \tilde \theta_0}{\tilde \theta_0^T \tilde \theta_0 +
      C_4}
    \\
    &+ \frac{C_3 (k - 1) \max\{\overline q, \overline \lambda\}^{2k} +
      C_2 \overline q^{2k}}{\inf_{\tilde \theta_0 \in (\real \ones_n)^\perp}
      (\tilde \theta_0^T \tilde \theta_0 + C_4)}
    \\
    &\hspace{-1in}= \overline \lambda^{2k} + C_3 C_4^{-1} (k - 1)
    \max\{\overline q, \overline \lambda\}^{2k} + C_2 C_4^{-1}
    \overline q^{2k}.
  \end{align*}
  By raising the right hand side of the above expression to the power
  $1/2k$ and taking the limit as $k \to \infty$, the
  constant/polynomial factors converge to 1 and the terms containing
  $\max\{\overline q, \overline \lambda\}$ dominate the sum, proving
  that $\mu \le \max\{\overline q, \overline \lambda\}$. Similarly, we
  can lower bound $\mu$ as follows. From~\eqref{eq:e-th2}, we have for
  all $k \in \integerpos$,
  \begin{align*}
    &\E\big[\hat \theta(k)^T \hat \theta(k)\big] \ge \theta_0^T
    \Ac^{2k} \theta_0
    \\
    &\Rightarrow \mu \ge \lim_{k \to \infty} \left(\sup_{\tilde
        \theta_0 \in (\real \ones_n)^\perp} \frac{\tilde \theta_0^T \Ac^{2k}
        \tilde \theta_0}{\tilde \theta_0^T \tilde \theta_0 + C_4}
    \right)^{1/2k} = \overline \lambda,
  \end{align*}
  and
  \begin{align*}
    &\E\big[\hat \theta(k)^T \hat \theta(k)\big] \ge \E[\eta(k)^T S
    \Pi_n^2 S \eta(k)] = \sum_{i = 1}^n C_{5i} q_i^{2k}
    \\
    &\Rightarrow \mu \ge \lim_{k \to \infty} \left(\sup_{\tilde
        \theta_0 \in (\real \ones_n)^\perp} \frac{\sum_{i = 1}^n C_{5i}
        q_i^{2k}}{\tilde \theta_0^T \tilde \theta_0 + C_4}
    \right)^{1/2k} = \overline q,
  \end{align*}
  where $C_{5i} = 2 c_i^2 (S \Pi_n^2 S)_{ii}$ for all $i \in
  \until{n}$. Therefore, $\mu \ge \max\{\overline q, \overline
  \lambda\}$, completing the proof. \qed
\end{pf}

Note that $\overline \lambda$ is the convergence rate of the
  noise-free (and non-private) Laplacian-based average consensus
  algorithm, while $\overline q$ is the worst-case decay rate of the
  noise sequence among the agents. From~\eqref{eq:mu}, the convergence
  rate $\mu$ is the larger of these two values, confirming our
  intuition that the slower rate among them is the bottleneck for
  convergence speed. Also, note that $\overline \lambda$ depends on
  the network topology $\Gc$ while $\overline q$ is independent of
  it.

\subsection{Accuracy and Differential
  Privacy}\label{sec:accuracy-privacy}

Having established the convergence properties of the
algorithm~\eqref{eq:dynamics}, here we characterize the extent to
which our design solves Problem~\ref{problem:main} by providing
guarantees on its accuracy and differential privacy.  The next result
elaborates on the statistical properties of the agreement value.

\begin{corollary}[Accuracy]\label{cor:accuracy}
  Under the hypotheses of Proposition~\ref{prop:convergence}, the
  convergence point $\theta_\infty$ is an unbiased estimate of
  $\Ave{\theta_0}$ with bounded dispersion,
  \begin{align}\label{eq:final-var}
    \text{var} \left \lbrace \theta_\infty \right \rbrace =
    \frac{2}{n^2} \sum_{i=1}^n \frac{s_i^2 c_i^2}{1 - q_i^2} .
  \end{align}
  As a result, the algorithm~\eqref{eq:dynamics}-\eqref{eq:noise} is $
  \Big( p, \tfrac{1}{n} \sqrt{\tfrac{2}{p} {\textstyle \sum_{i=1}^n}
    \tfrac{s_i^2 c_i^2}{1 - q_i^2}} \Big)$-accurate for any $p \in
  (0,1)$.
\end{corollary}
\begin{pf}
  Since noises are independent over time and among agents, we deduce
  from~\eqref{eq:dynamics-ave} that for any $k \in \integernonneg$, $
  E \lbrace \Ave{\theta(k)} \rbrace = \Ave{\theta_0}$ and
  \begin{align*}
    \text{var} \lbrace \Ave{\theta(k)} \rbrace = \frac{2}{n^2}
    \sum_{j=0}^k \sum_{i = 1}^n s_i^2 c_i^2 q_i^{2j} ,
  \end{align*}
  which establishes unbiasedness and bounded dispersion for any
  time. As $k \to \infty$, we get $ E \lbrace \theta_\infty \rbrace =
  \Ave{\theta_0}$ and
  \begin{align*}
    \text{var} \lbrace \theta_\infty \rbrace &= \frac{2}{n^2} \sum_{i
      = 1}^n \frac{s_i^2 c_i^2}{1 - q_i^2}.
  \end{align*}
  The $(p, r)$-accuracy follows directly by applying Chebyshev's
  inequality~\eqref{eq:Chebyshev} for $N = 1/\sqrt{p}$. \qed
\end{pf}

\begin{remark}\longthmtitle{Comparison with the Literature -- Cont'd} {\rm
    Proposition~\ref{prop:convergence} and Corollary~\ref{cor:accuracy}
    establish almost sure convergence, with the expected value of
    convergence being the average of the agents' initial states. In
    contrast, the results in~\citep{ZH-SM-GD:12} establish convergence
    in mean square, and the expected value of convergence depends on
    the network topology. In both cases, the accuracy radius $r$
    decreases with the number of agents as $O(1/\sqrt n)$.  }
  \oprocend
\end{remark}

The expression for $(p, r)$-accuracy in Corollary~\ref{cor:accuracy}
shows that one cannot obtain the ideal case of $(0, 0)$-accuracy, and
that $r$ is a decreasing function of $p$, with $r \to \infty$ as $p
\to 0$. This is an (undesirable) consequence of the lack of preservation of the
average under~\eqref{eq:dynamics} due to the term~$S \eta $. In turn,
the presence of this expression helps establish the differential
privacy of the algorithm with bounded, asymptotically vanishing noise,
as we show next.

\begin{proposition}[Differential Privacy]\label{prop:privacy}
  Under the hypotheses of Proposition~\ref{prop:convergence}, let
  \begin{align}\label{eq:epsilon}
  \epsilon_i = \delta \frac{q_i}{c_i(q_i - |s_i - 1|)} ,
  \end{align}
  for each $i \in \until{n}$, where $\delta$ is the adjacency bound
  in~\eqref{eq:adjacency}. Then, the algorithm preserves the
  $\epsilon_i$-differential privacy of agent $i$'s initial state for
  all $i \in \until{n}$. Consequently, the algorithm is
  $\epsilon$-differential private with $\epsilon = \max_i \epsilon_i$.
\end{proposition}
\begin{pf}
  Consider any pair of $\delta$-adjacent initial conditions
  $\theta_0^{(1)}$ and $\theta_0^{(2)}$ and an arbitrary set $\Ocal
  \subset (\real^n)^\integerpos$.  For any $k \in \integernonneg$, let
  \begin{align}\label{eq:regions}
    R_k^{(\ell)} = \lbrace \boldsymbol{\eta}_k \in \Omega_k \mid X_{k,
      \theta_0^{(\ell)}} (\boldsymbol{\eta}_k) \in \Ocal_k \rbrace,\
    \ell = 1,2 ,
  \end{align}
  where $\Omega_k = \real^{n(k+1)}$ is the sample space up to time
  $k$, $X_{k, \theta_0}$ is given in~\eqref{eq:mappings}, and $\Ocal_k
  \subseteq \real^{n(k + 1)}$ is the set composed by truncating the
  elements of $\Ocal$ to finite subsequences of length $k +
  1$. Then, by the continuity of probability~\cite[Theorem
    1.1.1.iv]{RD:10},
  \begin{multline}\label{eq:lim-int-noise}
    \pr \lbrace \boldsymbol{\eta} \in \Omega \mid X_{\theta_0^{(\ell)}}
    (\boldsymbol{\eta}) \in \Ocal \rbrace
    \\
    = \lim_{k \to \infty} {\int_{R_k^{(\ell)}}
      f_{n(k+1)}(\boldsymbol{\eta}_k^{(\ell)})d
      \boldsymbol{\eta}_k^{(\ell)}} , 
  \end{multline}
  for $ \ell = 1, 2$, where $f_{n(k+1)}$ is the $n(k+1)$-dimensional
  joint Laplace pdf given by
  \begin{align}\label{eq:nk-dim-pdf}
    f_{n(k+1)}(\boldsymbol{\eta}_k) = \prod_{i = 1}^n {\prod_{j = 0}^k \Lc
      (\eta_i(j) ; b_i(j))} .
  \end{align}
  Next, we define a bijection between $R_k^{(1)}$ and
  $R_k^{(2)}$. Without loss of generality, assume $\theta_{0,i_0}^{(2)}
  = \theta_{0,i_0}^{(1)} + \delta_1$ for some $i_0 \in \until{n}$,
  where $0 \le \delta_1 \le \delta$ and $\theta_{0,i}^{(2)} =
  \theta_{0,i}^{(1)}$ for all $i \neq i_0$. Then, for any
  $\boldsymbol{\eta}_k^{(1)} \in R_k^{(1)}$, define
  $\boldsymbol{\eta}_k^{(2)}$ by
  \begin{align*}
    \eta_i^{(2)}(j) &=
    \begin{cases}
      \eta_i^{(1)}(j) - (1 - s_i)^j \delta_1 , & \text{if $i =i_0$} ,
      \\
      \eta_i^{(1)}(j) , & \text{if $i \neq i_0$} ,
    \end{cases}
  \end{align*}
  for $j \in \lbrace 0 , \dots, k \rbrace$.  It is not difficult to see
  that $X_{k, \theta_0^{(1)}} (\boldsymbol{\eta}_k^{(1)}) = X_{k,
    \theta_0^{(2)}} (\boldsymbol{\eta}_k^{(2)})$, so
  $\boldsymbol{\eta}_k^{(2)} \in R_k^{(2)}$. Since the converse
  argument is also true, the above defines a bijection. Therefore, for
  any $\boldsymbol{\eta}_k^{(2)} \in R_k^{(2)}$ there exists a unique
  $(\boldsymbol{\eta}_k^{(1)}, \Delta \boldsymbol{\eta}_k) \in
  R_k^{(1)} \times \real^{n(k + 1)}$ such that
  \begin{align*}
    \boldsymbol{\eta}_k^{(2)} = \boldsymbol{\eta}_k^{(1)} + \Delta
    \boldsymbol{\eta}_k.
  \end{align*}
  Note that $\Delta \boldsymbol{\eta}_k$ is fixed and does not depend
  on $\boldsymbol{\eta}_k^{(2)}$.
  Thus, we can use a change of variables to get
  \begin{multline}\label{eq:lim-int-noise2}
    \pr \lbrace \boldsymbol{\eta} \in \Omega \mid X_{\theta_0^{(2)}}
    (\boldsymbol{\eta}) \in \Ocal \rbrace
    \\
    = \lim_{k \to \infty} {\int_{R_k^{(1)}}
      f_{n(k+1)}(\boldsymbol{\eta}_k^{(1)} + \Delta
      \boldsymbol{\eta}_k) d \boldsymbol{\eta}_k^{(1)}} .
  \end{multline}
  Comparing~\eqref{eq:lim-int-noise} for $\ell=1$
  with~\eqref{eq:lim-int-noise2}, we see that both integrals are over
  $R_k^{(1)}$ with different integrands. Dividing the integrands for
  any $\boldsymbol{\eta}_k^{(1)} \in R_k^{(1)}$ yields,
  \begin{align*}
    &\frac{f_{n(k+1)}
      (\boldsymbol{\eta}_k^{(1)})}{f_{n(k+1)}(\boldsymbol{\eta}_k^{(1)}
      + \Delta \boldsymbol{\eta}_k)}
    \\
    &\qquad = \frac{\prod_{i = 1}^n {\prod_{j=0}^k \Lc
        (\eta_i^{(1)}(j) ; b_i(j))}}{\prod_{i = 1}^n {\prod_{j=0}^k
        \Lc (\eta_i^{(1)}(j)+\Delta \eta_i(j) ; b_i(j))}}
    \\
    &\qquad = \frac {\prod_{j=0}^k \Lc (\eta_{i_0}^{(1)}(j) ;
      b_{i_0}(j))} {\prod_{j = 0}^k \Lc (\eta_{i_0}^{(1)}(j) + \Delta
      \eta_{i_0}(j) ; b_{i_0}(j))}
    \\
    &\qquad \le \prod_{j=0}^k e^{\frac{\left \lvert \Delta
          \eta_{i_0}(j) \right \rvert}{b_{i_0}(j)}} \le
    e^{\sum_{j=0}^k \frac{|1 - s_{i_0}|^j \delta}{c_{i_0} q_{i_0}^j}}
    \\
    &\Rightarrow f_{n(k+1)}(\boldsymbol{\eta}_k^{(1)}) \le
    e^{\frac{\delta}{c_{i_0}} \sum\limits_{j=0}^k \left( \frac{|1 -
          s_{i_0}|}{q_{i_0}} \right)^j}
    \!\!f_{n(k+1)}(\boldsymbol{\eta}_k^{(1)} + \Delta \boldsymbol{\eta}_k)
    .
  \end{align*}
  Due to~\eqref{eq:noise}, the geometric series in the exponent of the
  multiplicative term is convergent. Therefore, integrating both sides
  over $R_k^{(1)}$ and letting $k \to \infty$, we have
  \begin{multline*}
    \pr \lbrace \boldsymbol{\eta} \in \Omega \mid X_{\theta_0^{(1)}}
    (\boldsymbol{\eta}) \in \Ocal \rbrace
    \\
    \le e^{\delta \frac{q_{i_0}}{c_{i_0} (q_{i_0} - |1 - s_{i_0}|)}} \pr
    \lbrace \boldsymbol{\eta} \in \Omega \mid X_{\theta_0^{(2)}}
    (\boldsymbol{\eta}) \in \Ocal \rbrace ,
  \end{multline*}
  which establishes the $\epsilon_{i_0}$-differential privacy for
  agent~$i_0$. The fact the $i_0$ can be any agent
  establishes~\eqref{eq:epsilon}, while the last statement follows
  from Definition~\ref{def:privacy}. \qed
\end{pf}

Since the algorithm~\eqref{eq:dynamics}-\eqref{eq:noise} converges
almost surely (cf. Proposition~\ref{prop:convergence}) and is
differentially private (cf. Proposition~\ref{prop:privacy}),
Proposition~\ref{prop:bound} implies that it cannot achieve
$(0,0)$-accuracy, as noted above when discussing
Corollary~\ref{cor:accuracy}. The explicit privacy-accuracy
  trade-off is given by the relation between $\text{var}
  \{\theta_\infty\}$ and $\{\epsilon_i\}_{i = 1}^n$, i.e.,
  (c.f.~\eqref{eq:final-var},~\eqref{eq:epsilon})
\begin{align}\label{eq:priv-acc-to}
  \text{var} \left \lbrace \theta_\infty \right \rbrace = \frac{2
    \delta^2}{n^2} \sum_{i=1}^n \frac{s_i^2 q_i^2}{\epsilon_i^2 (q_i -
    |s_i - 1|)^2 (1 - q_i^2)},
\end{align}
so $\text{var} \left \lbrace \theta_\infty \right \rbrace$ increases
as any $\epsilon_i$ is decreased and vice versa. We optimize this
trade-off over $\{s_i, q_i\}_{i = 1}^n$ in
Section~\ref{section:optimal-noise-selection} and depict the optimal
trade-off curve for a test network in
Section~\ref{section:simulations}.

\begin{remark}[Laplacian Noise Distribution]
  {\rm Even though the choice of Laplacian noise in~\eqref{eq:noise}
    is not the only one that can be made to achieve differential
    privacy, it is predominant in the
    literature~\citep{CD-FM-KN-AS:06,CD:06}. The
    work~\citep{YW-ZH-SM-GED:14} shows that Laplacian noise is optimal
    (among all possible distributions) in the sense that it minimizes
    the entropy of the transmitted messages while preserving
    differential privacy.} \oprocend
\end{remark}

\begin{remark}\longthmtitle{Comparison with the Literature -- Cont'd}
  {\rm Proposition~\ref{prop:privacy} guarantees the
      $\epsilon_i$-differential privacy of agent $i$'s initial state
      independently of the noise levels chosen by other
      agents. Therefore, each agent can choose its own level of
    privacy, and even opt not to add any noise to its messages,
    without affecting the privacy of other agents.  In contrast,
    in~\citep{ZH-SM-GD:12}, agents need to agree on the level of
    privacy before executing the algorithm. In both cases, privacy is
    achieved against an adversary that can hear everything,
    independently of how it processes the information. In contrast,
    the algorithm in~\citep{YM-RMM:14,YM-RMM:15} assumes the adversary
    uses maximum likelihood estimation and only preserves the
    privacy of those agents who are sufficiently ``far'' from
    it in the graph (an agent is sufficiently far if the adversary
    cannot listen to it and all of its neighbors). The latter
      work uses a different notion of privacy based on the covariance
      of the maximum likelihood estimate which allows for guaranteed
      exact convergence, in the mean-square sense, to the true
      average.} \oprocend
\end{remark}

\subsection{Optimal Noise
  Selection}\label{section:optimal-noise-selection}
  
In this section, we discuss the effect on the algorithm's performance
of the free parameters present in our design. Given the trade-off
between accuracy and privacy, cf.~\eqref{eq:priv-acc-to},
we fix the privacy levels $\lbrace \epsilon_i \rbrace_{i = 1}^n$
constant and study the best achievable accuracy of the algorithm as a
function of the remaining free parameters.  Each agent $i \in
\until{n}$ gets to select the parameters $s_i$, $c_i$, $q_i$
determining the amount of noise introduced in the dynamics, with the
constraint that $ ( s_i, c_i, q_i) \in \Pc$, where
\begin{align*}
  \Pc = \lbrace ( s, c, q) \mid s \in (0, 2), c > 0, q \in (|s - 1|,
  1) \rbrace .
\end{align*}
Given the characterization of accuracy in
Corollary~\ref{cor:accuracy}, we consider as cost function the
variance of the agents' convergence point, i.e., $\theta_\infty$,
around $\Ave{\theta_0}$, giving
\begin{align}\label{eq:J}
  J (\lbrace s_i, c_i, q_i \rbrace_{i = 1}^n ) = \frac{2}{n^2}
  \sum_{i=1}^n \frac{s_i^2 c_i^2}{1 - q_i^2} .
\end{align}
The next result characterizes its global minimization.

\begin{proposition}\longthmtitle{Optimal Parameters for Variance
    Minimization}\label{prop:optimal}
  For the adjacency bound $\delta > 0$ and privacy levels $\lbrace
  \epsilon_i \rbrace_{i = 1}^n$ fixed, the optimal value of the variance of
  the agents' convergence point~is
  \begin{align*}
    J^* = \inf_{ \lbrace s_i, c_i, q_i \rbrace_{i = 1}^n \in \Pc^n} J
    (\lbrace s_i, c_i, q_i \rbrace_{i = 1}^n ) = 
    \frac{2 \delta^2}{n^2}
    \sum_{i = 1}^n \frac{1}{\epsilon_i^2}.
  \end{align*}
The infimum is not attained over $\Pc^n$ but approached as
  \begin{align}\label{eq:c_i}
    c_i = \delta \frac{q_i}{\epsilon_i(q_i - |s_i - 1|)} \quad , \quad
    s_i = 1,
  \end{align}
  and $q_i \to 0$ for all $i \in \until{n}$.
\end{proposition}
\begin{pf}
  For each $i \in \until{n}$, with the privacy level fixed, the
  expression~\eqref{eq:c_i} follows directly from~\eqref{eq:epsilon}.
  For convenience, we re-parameterize the noise decaying ratio~$q_i$
  as
  \begin{align}\label{eq:alpha}
    \alpha_i = \frac{q_i - \lvert s_i - 1 \rvert}{1 - \lvert s_i - 1
      \rvert} \in (0,1) .
  \end{align}
  Note that $ q_i = \alpha_i + (1 - \alpha_i) \lvert s_i - 1 \rvert$.
  Substituting~\eqref{eq:c_i} and~\eqref{eq:alpha} into~\eqref{eq:J},
  we obtain (with a slight abuse of notation, we also use $J$ to
  denote the resulting function),
  \begin{align*}
    &J (\lbrace s_i, \alpha_i \rbrace_{i = 1}^n ) = \frac{2}{n^2}
    \sum_{i = 1}^n \frac{\delta^2}{\epsilon_i^2} \phi(\alpha_i, s_i) ,
    \\
    &\phi(\alpha, s) = \frac{s^2 (\alpha + (1 - \alpha) \lvert s - 1
      \rvert )^2}{\alpha^2 (1 - |s - 1|)^2 \big[1 - (\alpha +
      (1 - \alpha) \lvert s - 1 \rvert )^2\big]}.
  \end{align*}
  Therefore, to minimize $J$, each agent has to independently minimize
  the same function $\phi$ of its local parameters $(\alpha_i, s_i)$
  over $D = (0, 1) \times (0, 2)$. Figure~\ref{fig:phi} illustrates
  the graph of this function over~$D$.
  \begin{figure}[htb]
    \includegraphics[width = 0.5\textwidth]{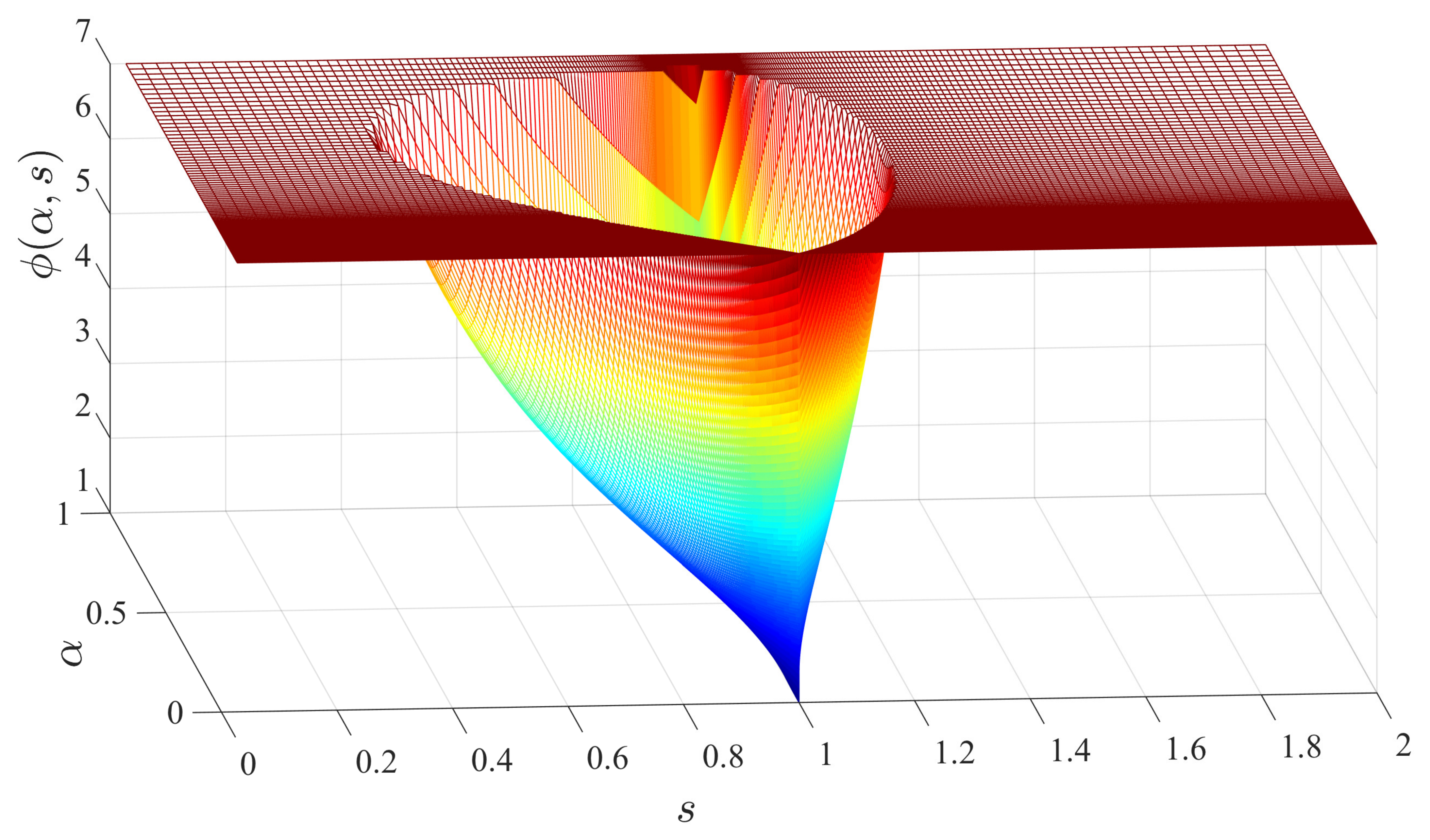}
    \caption{Local objective function~$\phi$ of each agent as a
      function of its parameters. $s$ is the noise-to-state gain and
      $\alpha$ is related to the noise decaying ratio. We cap
        the function values at $7$ for visualization purposes.  The
        function approaches its infimum as $\alpha \to 0$ while $s =
        1$.}\label{fig:phi}
  \end{figure}
  Since $D$ is not compact, the infimum might not be attained, and in
  fact, this is the case. It is easy verify that $\lim_{\alpha
      \to 0} \phi(\alpha, 1) = 1$.
  Now, for all $(\alpha, s) \in D$, $1 - (\alpha + (1 - \alpha) \lvert
  s - 1 \rvert)^2 < 1$ so
  \begin{align*}
    \phi(\alpha, s) > \phi_1^2(\alpha, s), \; \phi_1(\alpha, s) =
    \frac{(\alpha + (1 - \alpha) \lvert s - 1 \rvert)s}{\alpha
        (1 - |s - 1|)}.
  \end{align*}
  If $s \le 1$, then $\phi_1(\alpha, s) = s + \frac{1 -
      s}{\alpha} > 1$. If $s > 1$, then $\phi_1(\alpha, s) > 1 +
    \frac{s - 1}{\alpha(2 - s)} > 1$.
    Therefore, for all $(\alpha, s) \in D$, $\phi(\alpha, s) > 1$,
    which completes the proof. \qed
\end{pf}

Given that differential privacy is resilient to post-processing,
  an alternative design strategy to preserve the differential privacy
  of agents' initial states is to inject noise only at the initial
  time, $k=0$.  From~\eqref{eq:noise}, the introduction of a one-shot
  noise by agent~$i$ corresponds to $q_i = 0$ which is not feasible if
  $s_i \neq 1$. This can also be seen by rewriting~\eqref{eq:dynamics}
  as
  \begin{align*}
    \theta(k + 1) = (I - S)\theta(k) + (S - h L) x(k),
  \end{align*}
  so if $s_i \neq 1$ for any $i$,
  $\theta_i(k)$ directly (not only through $x_i(0)$) depend on
  $\theta_i(0)$. However, if $s_i = 1$,
  one can verify using a simplified version of the proof of
  Proposition~\ref{prop:privacy} that $q_i = 0$ also preserves
  $\epsilon_i$-differential privacy of $\theta_i(0)$ with $\epsilon_i
  = \frac{\delta}{c_i}$.  This results in a cost of
  \begin{align*}
    J = \frac{2}{n^2} \sum_{i = 1}^n c_i^2 = \frac{2 \delta^2}{n^2}
    \sum_{i = 1}^n \frac{1}{\epsilon_i^2} = J^*,
  \end{align*}
  showing that the optimal accuracy is also achieved by one-shot
  perturbation of the initial state at time $k = 0$ and injection of
  no noise thereafter.  A similar conclusion (that one-shot Laplace
  perturbation minimizes the output entropy) can be drawn
  from~\citep{YW-ZH-SM-GED:14}, albeit this is not explicitly
  mentioned therein.

  \begin{remark}\longthmtitle{Dynamic Average Consensus}
    {\rm
      In dynamic average
        consensus~\citep{HB-RAF-KML:10,MZ-SM:10-auto,SSK-JC-SM:15-ijrnc},
        agents seek to compute the average of individual exogenous,
        time-varying signals (the ``static'' average consensus
        considered here would be a special case corresponding to the
        exogenous signals being constant). In such scenarios, it is
        straightforward to show, using an argument similar to
        Proposition~\ref{prop:bound}, that one-shot perturbation would
        no longer preserve the differential privacy of time-varying
        input signals.
        The reason is that in this case, there is a recurrent flow of
        information at each node whose privacy can no longer be
        preserved with one-shot perturbation.  Sequential perturbation
        as in~\eqref{eq:message}-\eqref{eq:noise} is then necessary
        and the variance of the noise sequence has to dynamically
        depend on the rate of information flow to each node. Although
        the detailed design of such algorithms is beyond the scope of
        this work, such an algorithm can be designed following the
        idea of the sequential perturbation design of this work and
        the proof of its privacy in Proposition~\ref{prop:privacy}. To
        see this, note that (for $S \neq I_n$) we ``tune'' the amount
        of noise injection $\eta_i(k)$ so that the privacy of $(1 -
        s_i)^k \theta_{0, i}$ is preserved at each round $k \ge 1$,
        but $(1 - s_i)^k \theta_{0, i}$ is the amount of ``retained
        information'' of $\theta_{0, i}$ at round $k$ and plays the
        same role as $u(k)$ in the dynamic average consensus
        problem. \oprocend}
  \end{remark}

\section{Simulations}\label{section:simulations}

\begin{figure}[htb]
  \centering
  \includegraphics[width = 0.3\textwidth]{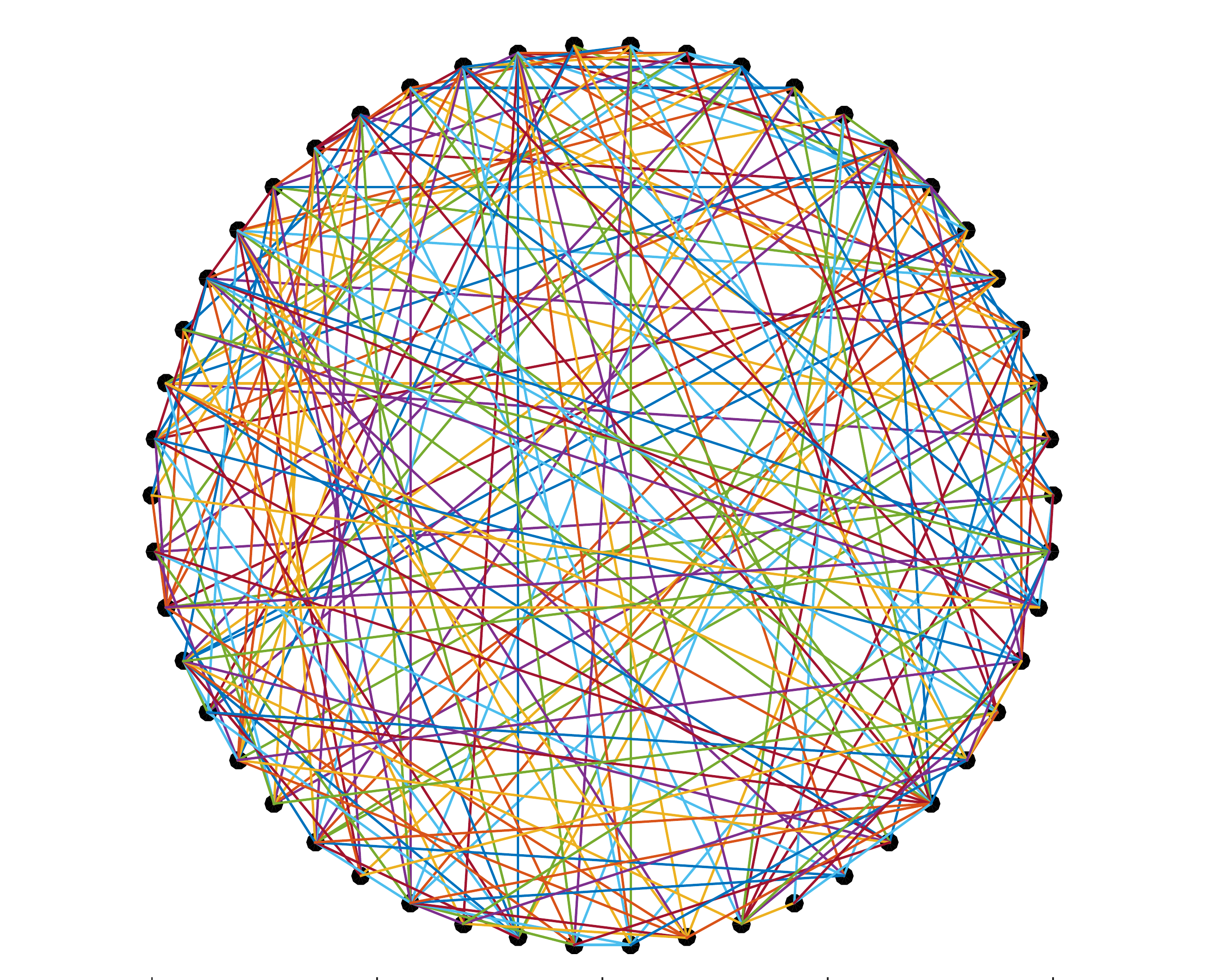}
  \caption{Random graph used for simulation.}\label{fig:graph}
\end{figure}

In this section, we report simulation results of the distributed
dynamics~\eqref{eq:dynamics}-\eqref{eq:noise} on a network of $n = 50$
agents. Figure~\ref{fig:graph} shows the random graph used
  throughout the section, where edge weights are i.i.d. and each one
  equals a sum of two i.i.d. Bernoulli random variables with $p =
  0.1$. The agents' initial states are also i.i.d. with distribution
  $\Nc(50, 100)$. As can be seen from~\eqref{eq:final-var}
and~\eqref{eq:epsilon}, neither accuracy nor privacy depend on the
initial values or the communication topology (albeit according
to~\eqref{eq:mu} the convergence rate depends on the latter). In
  all the simulations, $\delta = 1$ and $c_i = \delta q_i
  /\epsilon_i(q_i - |s_i - 1|)$ for all $i \in \until{n}$.

\begin{figure}[htb]
  \begin{minipage}{\linewidth}
    \subfloat[]{\includegraphics[width = 0.99\textwidth]{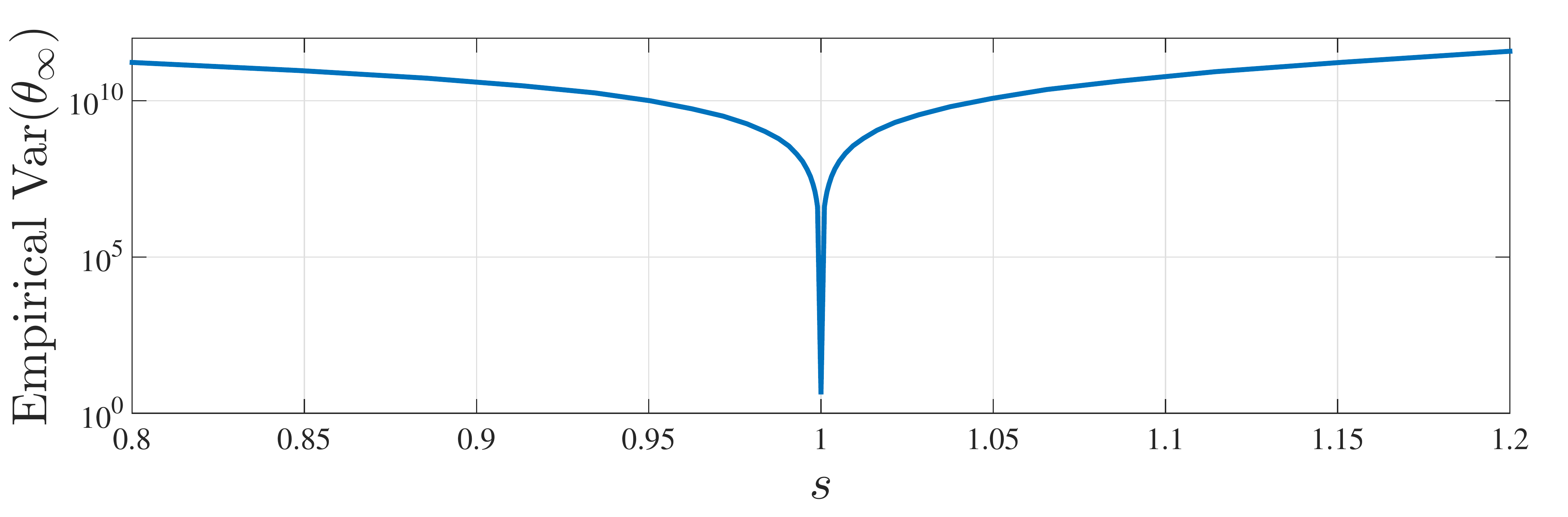}}
  \end{minipage}
  \begin{minipage}{\linewidth}
    \subfloat[]{\includegraphics[width = 0.99\textwidth]{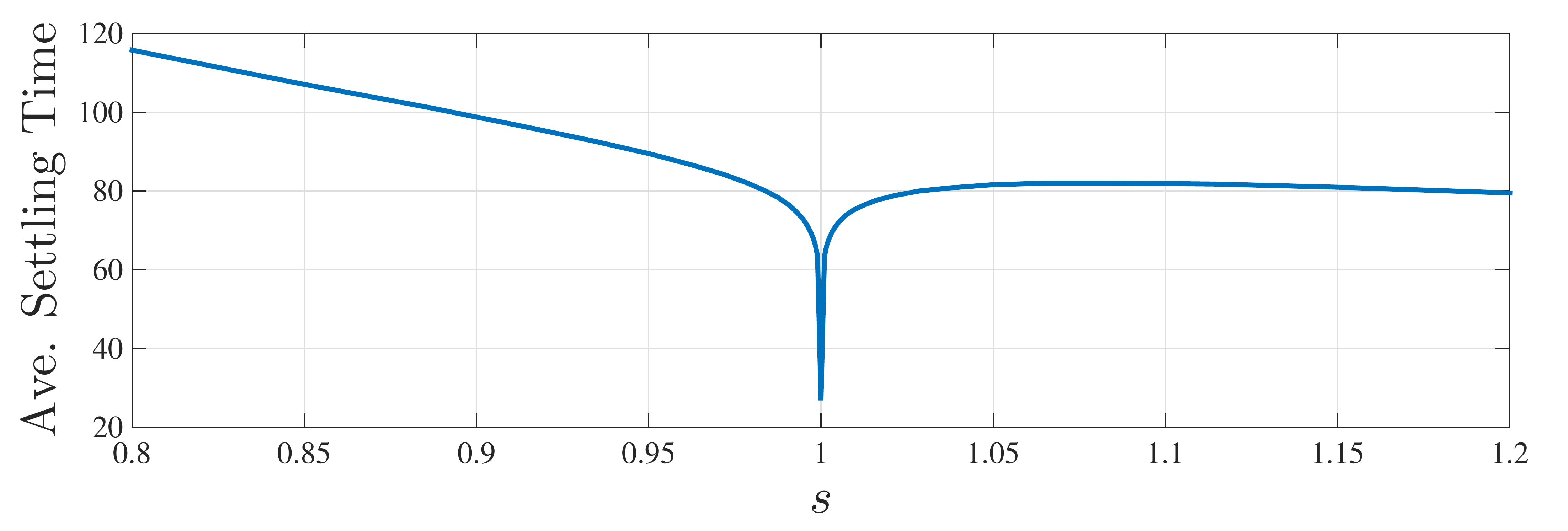}}
  \end{minipage}
  \caption{Executions of the
    algorithm~\eqref{eq:dynamics}-\eqref{eq:noise} for random topology
    and initial conditions.  (a) shows the empirical (i.e.,
      sample) variance of the convergence point and (b) shows the
    settling time. The trend in (a) validates
    Proposition~\ref{prop:optimal} while (b) shows the
      optimality of one-shot perturbation for convergence
      speed.}\label{fig:sweep}
\end{figure}

Figure~\ref{fig:sweep} depicts simulations with $\epsilon = 0.1
  \cdot \ones_n$ and $S = s I_n$ while sweeping~$s$ over $[0.8, 1.2]$
  with logarithmic step size. For each value of~$s$, we set $q_i =
  \alpha_i + (1 - \alpha_i) |s - 1|$ with $\alpha_i = 10^{-6}$ for
  each $i \in \until{n}$ and repeat the simulation $10^4$ times.  For
  each run, to capture the statistical properties of the convergence
  point, the graph topology and initial conditions are the same and
  only noise realizations change.  Figure~\ref{fig:sweep}(a) shows the
  empirical (sample) standard deviation of the convergence point as a
  function of~$s$, verifying the optimality of one-shot perturbation.
  In particular, notice the sensitivity of the accuracy to~$s$ close
  to $s = 1$.  Figure~\ref{fig:sweep}(b) shows the `settling time',
  defined as the number of rounds until convergence (measured by a
  tolerance of $10^{-2}$), as a function of~$s$.  The fastest
  convergence is achieved for $s = 1$, showing that one-shot noise is
  also optimal in the sense of convergence speed.  We have observed
  the same trends as in Figure~\ref{fig:sweep} for different random
  choices of initial conditions and network topologies.
  Note that the settling time depends on both the convergence rate and
  the initial distance from the convergence point $\|\theta(0) -
  \theta_\infty \ones_n\|$.  The former is constant at $\mu =
  \overline \lambda = 0.84$ for $s \in [0.8, 1.2]$.  The latter
  depends on $\{c_i\}_{i = 1}^n$, which in turn depend on $s$
  by~\eqref{eq:c_i}. This explains the trend observed in
  Figure~\ref{fig:sweep}(b).

\begin{figure}[htb]
\begin{minipage}{\linewidth}
    \subfloat[]{\includegraphics[width = 0.99\textwidth]{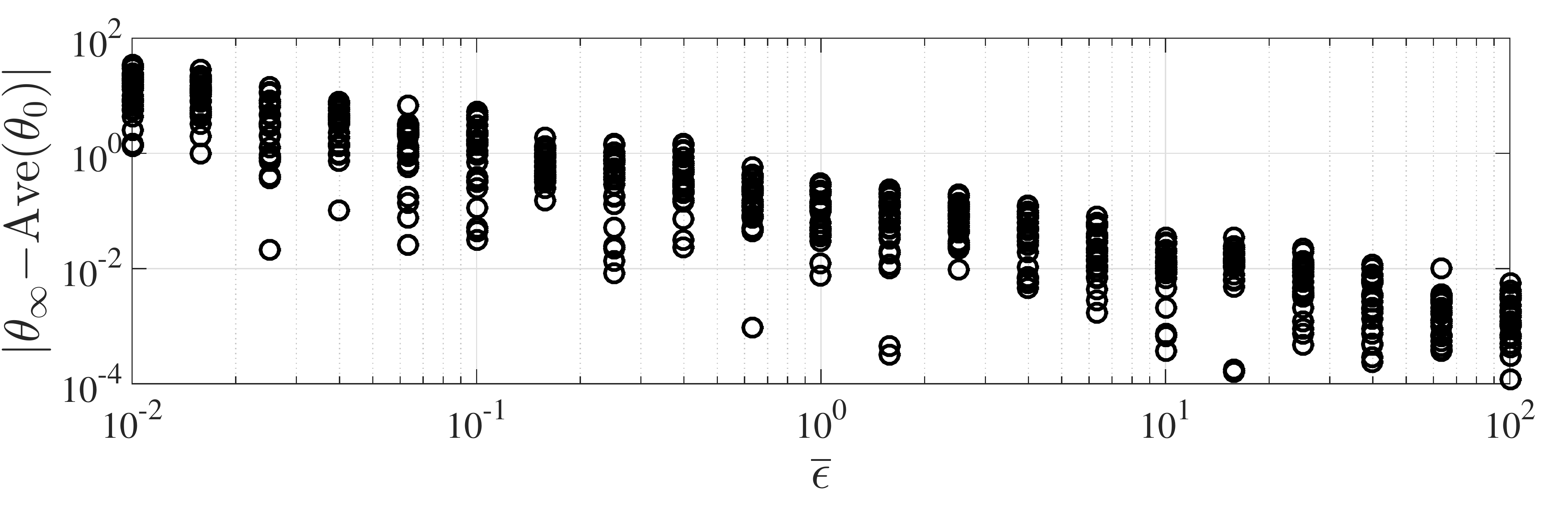}}
  \end{minipage}
  \begin{minipage}{\linewidth}
    \subfloat[]{\includegraphics[width = 0.99\textwidth]{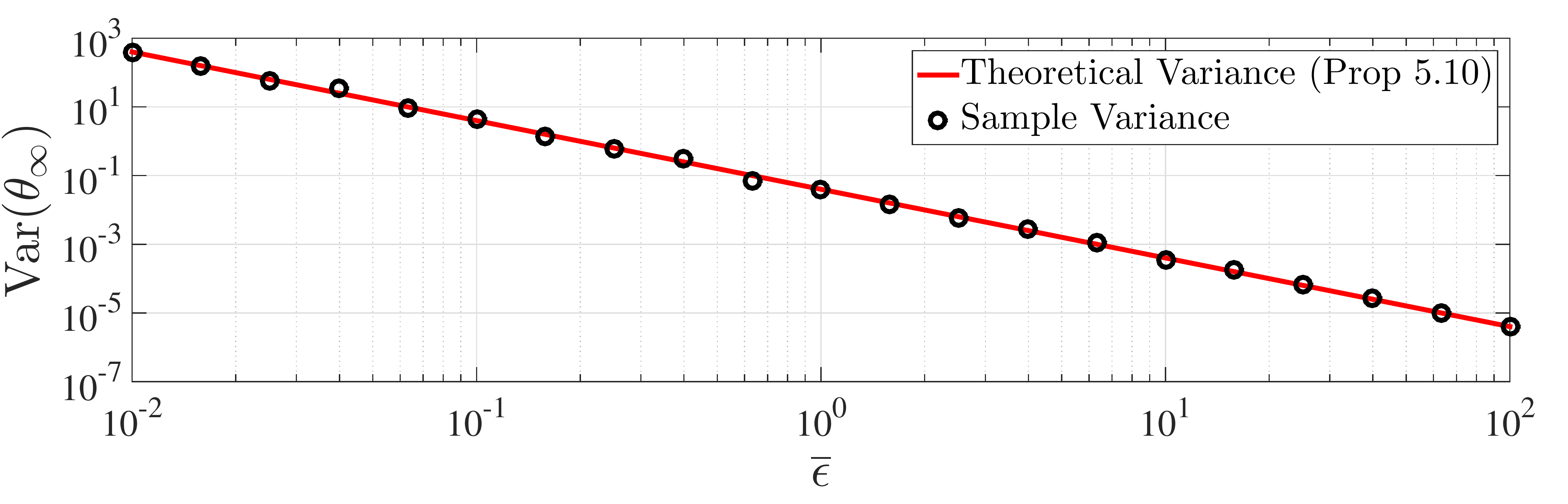}}
  \end{minipage}
  \caption{The privacy-accuracy trade-off for the proposed
    algorithm~\eqref{eq:dynamics}-\eqref{eq:noise} for random topology
    and initial conditions.  (a) shows the norm of the error for 25
    different realizations of the noise and (b) shows the sample
    variance over 100 noise realizations as well as the theoretical
    value provided by Proposition~\ref{prop:optimal}. The trend in
    both figures conforms with the theoretical characterization of
    $\theta_\infty$ given in
    Corollary~\ref{cor:accuracy}.}\label{fig:tradeoff}
\end{figure}

Figure~\ref{fig:tradeoff} depicts the privacy-accuracy trade-off for
the proposed algorithm. We have set $S = I_n$, $q = 0 \cdot \ones_n$,
and $\epsilon = \overline \epsilon \cdot \ones_n$ and then swept
$\overline \epsilon$ logarithmically over $[10^{-2}, 10^2]$. In
Figure~\ref{fig:tradeoff}(a), the algorithm is run 25 times for each
value of the $\overline \epsilon$ and the error $|\theta_\infty -
\Ave{\theta_0}|$ for each run is plotted as a circle. In
Figure~\ref{fig:tradeoff}(b), the sample variance of the convergence
point $\theta_\infty$ is shown as a function of $\overline \epsilon$
together with the theoretical value given in
Proposition~\ref{prop:optimal}. In both plots, we see an
inversely-proportional relationship between accuracy and privacy, as
expected.

\begin{figure}[htb]
  \includegraphics[width = 0.48\textwidth]{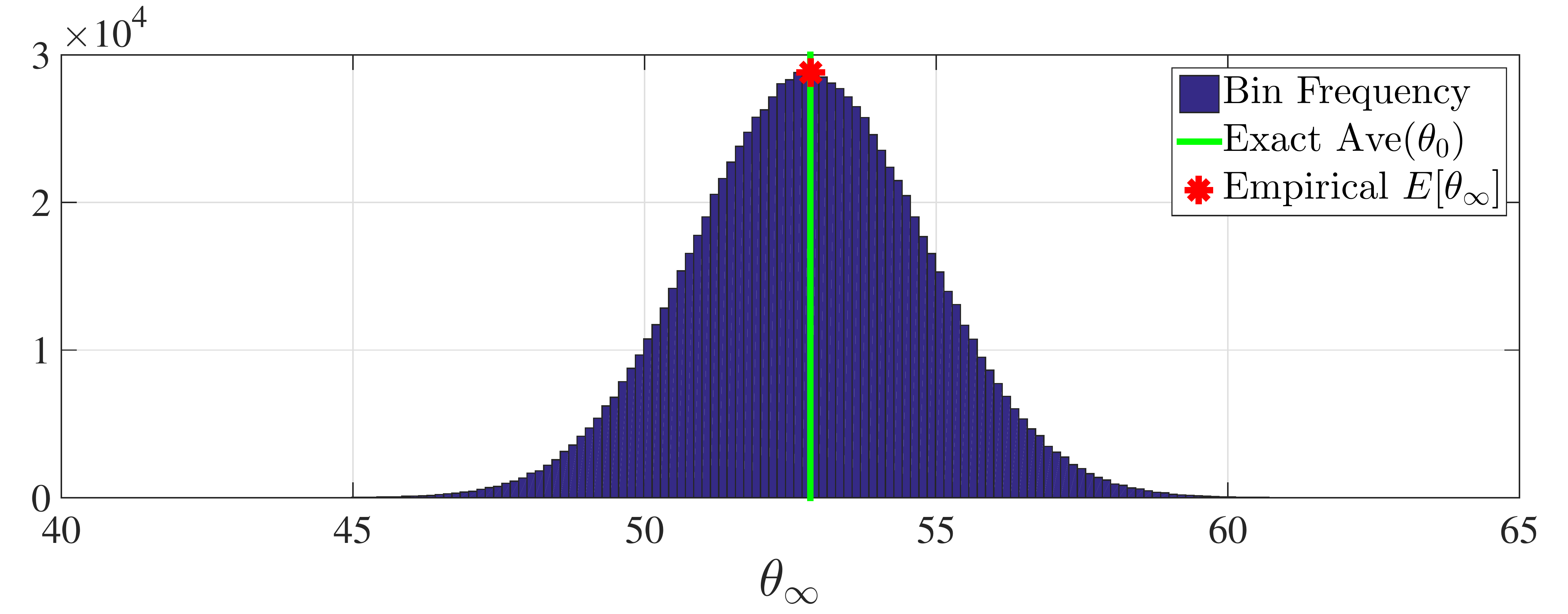}
  \caption{Statistical distribution of the convergence point. The
    sample mean (starred) matches the true average (green vertical
    line).}\label{fig:hist}
\end{figure}

Figure~\ref{fig:hist} shows the histogram of convergence points for
$10^6$ runs of the algorithm with $\epsilon = 0.1 \cdot \ones_n$, $S =
I_n$ and $q = 0 \cdot \ones_n$ (optimal accuracy). The distribution of
the convergence point is a bell-shaped curve with mean exactly at the
true average, in accordance with Corollary~\ref{cor:accuracy}.
Although the distribution of $\theta_\infty$ is provably non-Gaussian,
the central limit theorem, see e.g.,~\citep{RD:10}, implies that it is
very close to Gaussian since the number of agents is large.

\begin{figure}[htb]
  \includegraphics[width = 0.48\textwidth]{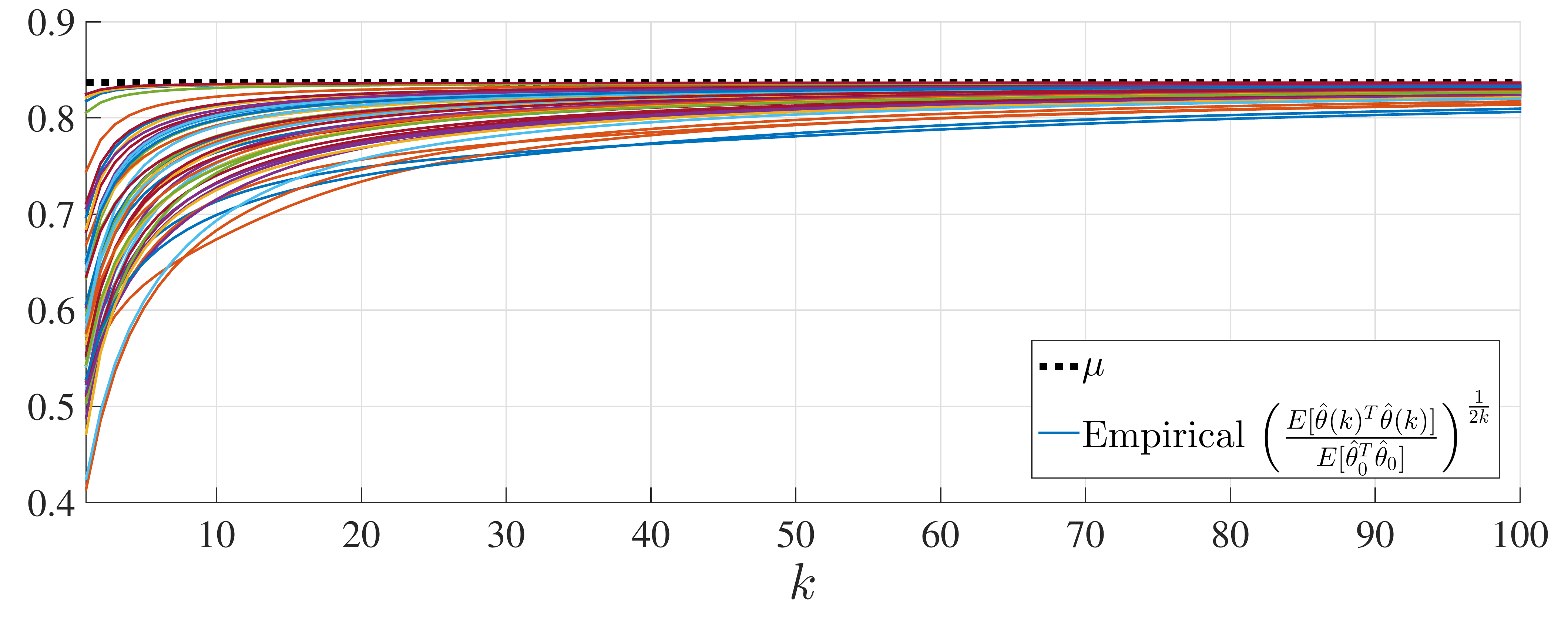}
  \caption{Illustration of the convergence rate of the
    algorithm~\eqref{eq:dynamics}-\eqref{eq:noise}. The limit of the
    supremum of the solid lines converges to the theoretical value of
    the exponential mean-square convergence rate $\mu$ given by
    Proposition~\ref{prop:conv-rate}. The curves with higher values
    correspond to initial states $\theta_0$ that are closer to the
    eigenvector of $I_n - h L - \Pi_n$ associated with $\overline
    \lambda$.}\label{fig:mu}
\end{figure}

Finally, Figure~\ref{fig:mu} illustrates the convergence rate of the
algorithm. Here, for $\epsilon = 0.1 \ones_n$, $S = 0.9 I_n$, $q = 0.2
  \ones_n$, and the same topology as in the
previous plots, the initial agents states are randomly selected and
the whole algorithm is run 100 times with different noise realizations
$\boldsymbol{\eta}$, each time until $100$ iterations.  For each value
of initial states and each $k \in \until{100}$, we empirically
approximate the quantity \vspace*{-4pt}
\begin{align*}
  \left( \frac{\E\big[(\theta(k) - \theta_\infty \ones_n)^T (\theta(k)
      - \theta_\infty \ones_n) \big]}{\E\big[(\theta(0) -
      \theta_\infty \ones_n)^T (\theta(0) - \theta_\infty
      \ones_n)\big]} \right)^{1/2k}
\end{align*}
by taking the sample mean instead of the expectation in the numerator
and denominator. We repeat this whole process 50 times for different
random initial conditions and plot the result, together with the
theoretical value of $\mu$ (which in this case equals $\overline
\lambda$) given by Proposition~\ref{prop:conv-rate}. As
Figure~\ref{fig:mu} shows, the supremum of the resulting curves
converges to $\mu$ as $k \to \infty$, as expected.

\section{Conclusions}\label{section:conclusions}

We have studied the problem of multi-agent average consensus subject
to the differential privacy of agents' initial states.  We have showed
that the requirement of differential privacy cannot be satisfied if
agents' states weakly converge to the exact average of their initial
states. This result suggests that the most one can expect of a
differentially private consensus algorithm is that the consensus value
is unbiased, i.e., its expected value is the true average, and the
variance is minimized. We have designed a linear consensus algorithm
that meets this objective, and have carefully characterized its
convergence, accuracy, and differential privacy properties.  Future
work will include the investigation of the limitations and advantages
of differential privacy in multi-agent systems, the extension of the
results to dynamic average consensus, distributed optimization,
filtering, and estimation, and the design of algorithms for privacy
preservation of the network structure and other parameters such as
edge weights and vertex degrees.

\section*{Acknowledgments}
This work was supported by NSF Award CNS-1329619.

{
\small
\bibliographystyle{plainnat}

\begin{thebibliography}{31}
\providecommand{\natexlab}[1]{#1}
\providecommand{\url}[1]{\texttt{#1}}
\expandafter\ifx\csname urlstyle\endcsname\relax
  \providecommand{\doi}[1]{doi: #1}\else
  \providecommand{\doi}{doi: \begingroup \urlstyle{rm}\Url}\fi

\bibitem[Bai et~al.(2010)Bai, Freeman, and Lynch]{HB-RAF-KML:10}
H.~Bai, R.~A. Freeman, and K.~M. Lynch.
\newblock Robust dynamic average consensus of time-varying inputs.
\newblock In \emph{{IEEE} Conf.\ on Decision and Control}, pages 3104--3109,
  Atlanta, GA, December 2010.

\bibitem[Bullo et~al.(2009)Bullo, Cort{\'e}s, and
  Mart{\'\i}nez]{FB-JC-SM:08cor}
F.~Bullo, J.~Cort{\'e}s, and S.~Mart{\'\i}nez.
\newblock \emph{Distributed Control of Robotic Networks}.
\newblock Applied Mathematics Series. Princeton University Press, 2009.
\newblock ISBN 978-0-691-14195-4.
\newblock Electronically available at \url{http://coordinationbook.info}.

\bibitem[Duan et~al.(2015)Duan, He, Cheng, Mo, and Chen]{XD-etal:15}
X.~Duan, J.~He, P.~Cheng, Y.~Mo, and J.~Chen.
\newblock Privacy preserving maximum consensus.
\newblock In \emph{{IEEE} Conf.\ on Decision and Control}, pages 4517--4522,
  Osaka, 2015.

\bibitem[Durrett(2010)]{RD:10}
R.~Durrett.
\newblock \emph{Probability: Theory and Examples}.
\newblock Series in Statistical and Probabilistic Mathematics. Cambridge
  University Press, 4th edition, 2010.
\newblock ISBN 9780521765398.

\bibitem[Dwork(2006)]{CD:06}
C.~Dwork.
\newblock Differential privacy.
\newblock In \emph{Proceedings of the 33rd International Colloquium on
  Automata, Languages and Programming (ICALP)}, pages 1--12, Venice, Italy,
  July 2006.

\bibitem[Dwork and Roth(2014)]{CD-AR:14}
C.~Dwork and A.~Roth.
\newblock The algorithmic foundations of differential privacy.
\newblock \emph{Found. Trends Theor. Comput. Sci.}, 9\penalty0 (3-4):\penalty0
  211--407, August 2014.

\bibitem[Dwork et~al.(2006)Dwork, McSherry, Nissim, and Smith]{CD-FM-KN-AS:06}
C.~Dwork, F.~McSherry, K.~Nissim, and A.~Smith.
\newblock Calibrating noise to sensitivity in private data analysis.
\newblock In \emph{Proceedings of the 3rd Theory of Cryptography Conference},
  pages 265--284, New York, NY, March 2006.

\bibitem[G{\"u}nd{\"u}z et~al.(2010)G{\"u}nd{\"u}z, Erkip, and
  Poor]{DG-EE-HVP:10}
D.~G{\"u}nd{\"u}z, E.~Erkip, and H.V. Poor.
\newblock Source coding under secrecy constraints.
\newblock In \emph{Securing Wireless Communications at the Physical Layer},
  pages 173--199. Springer US, Boston, MA, 2010.

\bibitem[Han et~al.(2014)Han, Topcu, and Pappas]{SH-UT-GJP:14}
S.~Han, U.~Topcu, and G.~J. Pappas.
\newblock Differentially private convex optimization with piecewise affine
  objectives.
\newblock In \emph{{IEEE} Conf.\ on Decision and Control}, pages 2160--2166,
  Los Angeles, CA, December 2014.

\bibitem[Huang et~al.(2012)Huang, Mitra, and Dullerud]{ZH-SM-GD:12}
Z.~Huang, S.~Mitra, and G.~Dullerud.
\newblock Differentially private iterative synchronous consensus.
\newblock In \emph{Proceedings of the 2012 ACM workshop on Privacy in the
  electronic society}, pages 81--90, New York, NY, 2012.

\bibitem[Huang et~al.(2014)Huang, Wang, Mitra, and Dullerud]{ZH-YW-SM-GED:14}
Z.~Huang, Y.~Wang, S.~Mitra, and G.~E. Dullerud.
\newblock On the cost of differential privacy in distributed control systems.
\newblock In \emph{Proceedings of the 3rd International Conference on High
  Confidence Networked Systems (HiCoNS)}, pages 105--114, Berlin, Germany,
  April 2014.

\bibitem[Huang et~al.(2015)Huang, Mitra, and Vaidya]{ZH-SM-NV:15}
Z.~Huang, S.~Mitra, and N.~Vaidya.
\newblock Differentially private distributed optimization.
\newblock In \emph{Proceedings of the 2015 International Conference on
  Distributed Computing and Networking}, Pilani, India, January 2015.

\bibitem[Jiang and Wang(2001)]{ZPJ-YW:01}
Z.-P. Jiang and Y.~Wang.
\newblock Input-to-state stability for discrete-time nonlinear systems.
\newblock \emph{Automatica}, 37\penalty0 (6):\penalty0 857--869, 2001.

\bibitem[Kairouz et~al.(2015)Kairouz, Oh, and Viswanath]{PK-SO-PV:15}
P.~Kairouz, S.~Oh, and P.~Viswanath.
\newblock Secure multi-party differential privacy.
\newblock In \emph{Advances in Neural Information Processing Systems 28}, pages
  2008--2016. Curran Associates, Inc., 2015.

\bibitem[Kefayati et~al.(2007)Kefayati, Talebi, Khalaj, and
  Rabiee]{MK-MST-BHK-HRR:07}
M.~Kefayati, M.~S. Talebi, B.~H. Khalaj, and H.~R. Rabiee.
\newblock Secure consensus averaging in sensor networks using random offsets.
\newblock In \emph{IEEE Intern. Conf. on Telec., and Malaysia Intern. Conf. on
  Communications}, pages 556--560, Penang, May 2007.

\bibitem[Kia et~al.(2015)Kia, Cort{\'e}s, and
  Mart{\'\i}nez]{SSK-JC-SM:15-ijrnc}
S.~S. Kia, J.~Cort{\'e}s, and S.~Mart{\'\i}nez.
\newblock Dynamic average consensus under limited control authority and privacy
  requirements.
\newblock \emph{International Journal on Robust and Nonlinear Control},
  25\penalty0 (13):\penalty0 1941--1966, 2015.

\bibitem[Manitara and Hadjicostis(2013)]{NEM-CNH:13}
N.~E. Manitara and C.~N. Hadjicostis.
\newblock Privacy-preserving asymptotic average consensus.
\newblock In \emph{{E}uropean {C}ontrol {C}onference}, pages 760--765, Zurich,
  Switzerland, 2013.

\bibitem[Mesbahi and Egerstedt(2010)]{MM-ME:10}
M.~Mesbahi and M.~Egerstedt.
\newblock \emph{Graph Theoretic Methods in Multiagent Networks}.
\newblock Applied Mathematics Series. Princeton University Press, 2010.

\bibitem[Mo and Murray(2014)]{YM-RMM:14}
Y.~Mo and R.~M. Murray.
\newblock Privacy preserving average consensus.
\newblock In \emph{{IEEE} Conf.\ on Decision and Control}, pages 2154--2159,
  Los Angeles, CA, December 2014.

\bibitem[Mo and Murray(2015)]{YM-RMM:15}
Y.~Mo and R.~M. Murray.
\newblock Privacy preserving average consensus.
\newblock \emph{IEEE Transactions on Automatic Control}, 2015.
\newblock Submitted, available at
  \url{http://yilinmo.github.io/public/papers/tac2014privacy.pdf}.

\bibitem[Mukherjee et~al.(2014)Mukherjee, Fakoorian, Huang, and
  Swindlehurst]{AM-SAAF-JH-ALS:14}
A.~Mukherjee, S.A.A. Fakoorian, J.~Huang, and A.L. Swindlehurst.
\newblock Principles of physical layer security in multiuser wireless networks:
  A survey.
\newblock \emph{IEEE Communications Surveys \& Tutorials}, 16\penalty0
  (3):\penalty0 1550--1573, 2014.

\bibitem[Nozari et~al.(2015)Nozari, Tallapragada, and
  Cort{\'e}s]{EN-PT-JC:15-necsys}
E.~Nozari, P.~Tallapragada, and J.~Cort{\'e}s.
\newblock Differentially private average consensus with optimal noise
  selection.
\newblock \emph{IFAC-PapersOnLine}, 48\penalty0 (22):\penalty0 203--208, 2015.
\newblock {\it IFAC Workshop on Distributed Estimation and Control in Networked
  Systems}, Philadelphia, PA.

\bibitem[Nozari et~al.(2017)Nozari, Tallapragada, and
  Cort{\'e}s]{EN-PT-JC:17-tcns}
E.~Nozari, P.~Tallapragada, and J.~Cort{\'e}s.
\newblock Differentially private distributed convex optimization via functional
  perturbation.
\newblock \emph{IEEE Transactions on Control of Network Systems}, 2017.
\newblock To appear.

\bibitem[Ny and Pappas(2014)]{JLN-GJP:14}
J.~L. Ny and G.~J. Pappas.
\newblock Differentially private filtering.
\newblock \emph{IEEE Transactions on Automatic Control}, 59\penalty0
  (2):\penalty0 341--354, 2014.

\bibitem[Olfati-Saber et~al.(2007)Olfati-Saber, Fax, and
  Murray]{ROS-JAF-RMM:07}
R.~Olfati-Saber, J.~A. Fax, and R.~M. Murray.
\newblock Consensus and cooperation in networked multi-agent systems.
\newblock \emph{Proceedings of the IEEE}, 95\penalty0 (1):\penalty0 215--233,
  2007.

\bibitem[Papoulis and Pillai(2002)]{AP-SUP:02}
A.~Papoulis and S.~U. Pillai, editors.
\newblock \emph{Probability, Random Variables and Stochastic Processes}.
\newblock McGraw-Hill, 2002.
\newblock ISBN 0073660116.

\bibitem[Pettai and Laud(2015)]{MP-PL:15}
M.~Pettai and P.~Laud.
\newblock Combining differential privacy and secure multiparty computation.
\newblock In \emph{Proceedings of the 31st Annual Computer Security
  Applications Conference}, ACSAC 2015, pages 421--430. ACM, 2015.

\bibitem[Ren and Beard(2008)]{WR-RWB:08}
W.~Ren and R.~W. Beard.
\newblock \emph{Distributed Consensus in Multi-Vehicle Cooperative Control}.
\newblock Communications and Control Engineering. Springer, 2008.
\newblock ISBN 978-1-84800-014-8.

\bibitem[Tanaka and Sandberg(2015)]{TT-HS:15}
T.~Tanaka and H.~Sandberg.
\newblock {SDP}-based joint sensor and controller design for
  information-regularized optimal {LQG} control.
\newblock In \emph{{IEEE} Conf.\ on Decision and Control}, pages 4486--4491,
  Osaka, 2015.

\bibitem[Wang et~al.(2014)Wang, Huang, Mitra, and Dullerud]{YW-ZH-SM-GED:14}
Y.~Wang, Z.~Huang, S.~Mitra, and G.~E. Dullerud.
\newblock Entropy-minimizing mechanism for differential privacy of
  discrete-time linear feedback systems.
\newblock In \emph{{IEEE} Conf.\ on Decision and Control}, pages 2130--2135,
  Los Angeles, CA, December 2014.

\bibitem[Zhu and Mart{\'\i}nez(2010)]{MZ-SM:10-auto}
M.~Zhu and S.~Mart{\'\i}nez.
\newblock Discrete-time dynamic average consensus.
\newblock \emph{Automatica}, 46\penalty0 (2):\penalty0 322--329, 2010.

\end{thebibliography}

}

\end{document}